\documentclass[a4paper,11pt,leqno]{article}

\usepackage[dvips]{graphics}
\usepackage[matrix,arrow,tips,curve]{xy}
\usepackage[english]{babel}
\usepackage{amsmath}
\usepackage{amssymb,amsthm}
\usepackage{enumerate}
\usepackage{psfrag}
 \psfrag{X}{\fontsize{20}{20}$X$}
 \psfrag{Y2}{\fontsize{20}{20}$Y_2$}
 \psfrag{Y1}{\fontsize{20}{20}$Y_1$}
 \psfrag{E1}{\fontsize{20}{20}$E_1$}
 \psfrag{E2}{\fontsize{20}{20}$E_2$}
 \psfrag{E3}{\fontsize{20}{20}$\widetilde{E}_2$}
 \psfrag{E4}{\fontsize{20}{20}$E_4$}
 \psfrag{E6}{\fontsize{20}{20}$D'$}
 \psfrag{ph2}{\fontsize{20}{20}$\ph_2$}
 \psfrag{ph1}{\fontsize{20}{20}$\ph_1$}
 \psfrag{t1}{\fontsize{20}{20}$\ph_1(E_1)$}
 \psfrag{t2}{\fontsize{20}{20}$\ph_1(E_2)$}
 \psfrag{D}{\fontsize{20}{20}$D$}
 \psfrag{A}{\fontsize{20}{20}$A$}
 \psfrag{P}{\fontsize{20}{20}$\pr^2$}
 \psfrag{R3}{\fontsize{20}{20}$R_3$}
 \psfrag{R2}{\fontsize{20}{20}$R_2$}

\newtheorem{thm}{Theorem}[section]
\newtheorem*{thm*}{Theorem}
\newtheorem{lemma}[thm]{Lemma}
\newtheorem{proposition}[thm]{Proposition}
\newtheorem{corollary}[thm]{Corollary}

\theoremstyle{definition}

\newtheorem{remark}[thm]{Remark}
\newtheorem{example}[thm]{Example}

\newcommand{\la}{\longrightarrow}

\newcommand{\ma}{\mathcal}
\newcommand{\ph}{\varphi}
\newcommand{\pr}{\mathbb{P}}

\newcommand{\Z}{\mathbb{Z}}
\newcommand{\Q}{\mathbb{Q}}
\newcommand{\R}{\mathbb{R}}

\newcommand{\N}{\ma{N}_1}

\newcommand{\Sing}{\operatorname{Sing}}

\newcommand{\NE}{\operatorname{NE}}

\newcommand{\Exc}{\operatorname{Exc}}

\newcommand{\Lo}{\operatorname{Locus}}

\newlength{\Mheight}
\newlength{\cwidth}
\title{On Fano manifolds with a birational contraction sending a
  divisor to a curve}
\author{C.\ 
Casagrande\footnote{Dipartimento di Matematica, Universit\`a di Pavia,
via Ferrata 1, 
 27100 Pavia - Italy,
 cinzia.casagrande@unipv.it}}
\date{July 15, 2008}

\begin{document}
\maketitle
\section{Introduction}
\noindent Let $X$ be a smooth, complex Fano variety of dimension $n\geq
4$. The Picard number $\rho_X$ of $X$ is equal to the second Betti
number of $X$, and is bounded in any fixed dimension; 
however the maximal value is unknown even in dimension $4$.

Bounds on $\rho_X$ are known when $X$ has some special extremal
contraction. For instance if $X$ has a birational elementary
contraction sending a divisor to a point, then
$\rho_X\leq 3$ (\cite[Prop.\ 5]{toru}, see also Prop.\ \ref{1}). 
In fact such $X$ are classified in
the toric case \cite{bonavero3}, in the case of a blow-up of a point
\cite{bonwisncamp}, and more generally when the exceptional divisor is
$\pr^{n-1}$ \cite{toru}.

Concerning the fiber type case, we know that $\rho_X\leq 11$ when $X$
has an elementary contraction onto a surface or a threefold
\cite[Th.\ 1.1]{fanos}. 

Here we consider the case of a 
birational elementary contraction 
of type $(n-1,1)$, that is,
sending a
divisor to a curve. Such Fano varieties have been classified in the toric case
by H.\ Sato
 \cite{sato2}, and T.\ Tsukioka has obtained
classification results 
for some cases 
\cite{toruCRAS,toru4fold} (see Rem.~\ref{tt}).
Our main result is the following.
\begin{thm}\label{5}
Let $X$ be a smooth Fano variety of dimension $n\geq 4$, 
and suppose that $X$ has a birational elementary contraction 
sending a
divisor $E$ to a curve. 

Then
$\rho_X\leq 5$.
Moreover if $\rho_X=5$ we have $E\cong W\times\pr^1$,
 $W$ a smooth Fano variety, 
and there exist: 
\begin{enumerate}[$\bullet$]
\item a smooth projective variety $Y$ with $\rho_Y=4$,
such that
$X$ is the blow-up of $Y$ in a subvariety isomorphic to $W$, with
exceptional divisor $E$;
\item  a smooth Fano variety $Z$
with $\rho_Z=3$, having a birational elementary contraction 
sending a
divisor $E_Z$ to a curve,
 such that $X$ is the blow-up of $Z$ in two fibers of such
 contraction, and $E$ is the proper transform of $E_Z$.
\end{enumerate}
\end{thm}
\noindent 
This Theorem follows from Th.~\ref{positive} and
Prop.~\ref{11/4}. There are examples with $\rho_X=5$ in every
dimension $n\geq 4$, see Ex.~\ref{ultimo}.

We give some applications to the $4$-dimensional case.
\begin{corollary}\label{dim4}
Let $X$ be a smooth Fano $4$-fold, 
and suppose that $X$ has a birational elementary contraction sending a
divisor $E$ to a curve. 

 Then
$\rho_X\leq 5$, and if $\rho_X=5$ we have one of the possibilities:
\begin{enumerate}[$(i)$]
\item $E\cong\pr^2\times\pr^1$,
  $\mathcal{N}_{E/X}\cong\mathcal{O}(-1,-1)$;
\item $E\cong\pr^2\times\pr^1$,
  $\mathcal{N}_{E/X}\cong\mathcal{O}(-2,-1)$;
\item $E\cong\pr^1\times\pr^1\times\pr^1$,
  $\mathcal{N}_{E/X}\cong\mathcal{O}(-1,-1,-1)$, and two of the
  rulings are numerically equivalent in $X$.
\end{enumerate}
\end{corollary}
\begin{corollary}\label{4folds}
Let $X$ be a smooth Fano $4$-fold. Then one of the
following holds:
\begin{enumerate}[$(i)$]
\item $\rho_X\leq 6$;
\item $X$ is a product and $\rho_X\leq 11$;
\item
every
elementary contraction of $X$ is birational of type $(3,2)$ or $(2,0)$.
\end{enumerate}
\end{corollary}

\medskip

We explain the technique used
 to prove Th.~\ref{5}. 
Given the
divisor $E$, the classical approach is to
 choose an extremal ray $R$ of $\NE(X)$ such that
$E\cdot R>0$, and study the associated contraction.
Anyway this is not enough to get a bound
on $\rho_X$ in all cases, in particular when $R$ is small.
One has to iterate this procedure and
run a ``Mori
program'' for $-E$, that is, to contract or flip birational extremal
rays having positive intersection with $E$, until one gets a fiber type
contraction. This is possible thanks to \cite{BCHM}, where it is shown
that Fano varieties are Mori dream spaces, and \cite{hukeel}, where
properties of Mori dream spaces are studied. 

In section
\ref{D} we use this method
to study a Fano variety $X$ containing  
a prime divisor $D$ such that the numerical classes of
curves contained in $D$ span a $2$-dimensional linear subspace in
$\N(X)$. This is enough to get $\rho_X\leq 3$ in some cases
(see Th.~\ref{nuovo}).

Then in section \ref{31} we consider the exceptional divisor
 $E$ of an
elementary contraction $\ph$ of type $(n-1,1)$.
We apply to $E$ the results of the preceding section, and
we need a detailed analysis of the geometry of $E$ and $X$
to conclude. We first show that if there is a unique
extremal ray having negative intersection with $E$ (corresponding to
$\ph$), then $\rho_X\leq
4$ (Th.~\ref{positive}). Then we consider  the case where there is
a second extremal ray $R$ such that $E\cdot R<0$, and show that
$\rho_X\leq 5$ (Prop.~\ref{11/4}).
Finally we give some examples with $\rho_X=5$.
\section{Preliminaries}\label{prel}
In this section we recall
some notions and results that we need in the sequel. 

\medskip

\noindent {\bf Contractions.}
Let $X$ be a normal irreducible variety of dimension $n$.
A \emph{contraction} of $X$ is a projective morphism $\ph\colon X\to
Y$, with connected fibers, onto a normal variety $Y$
(without hypotheses on the anticanonical degree of curves in
fibers). 
We say that $\ph$
is \emph{of type $(a,b)$} if $\dim \Exc(\ph)=a$ and $\dim
\ph(\Exc(\ph))=b$, where $\Exc(\ph)$ is the exceptional locus of $\ph$. 

Suppose that $X$ has terminal singularities, so that $K_X$ is
$\Q$-Cartier.
A contraction $\ph$ 
is a \emph{Mori contraction} if $-K_X$ is $\ph$-ample.

\medskip

\noindent {\bf Numerical equivalence classes and the cone of curves.}
Let $X$ be an irreducible projective variety.
 We denote by $\N(X)$ the vector
space of 1-cycles in $X$, with real coefficients, modulo numerical
equivalence. Its dimension is the Picard number $\rho_X$ of $X$. 
The cone of curves $\NE(X)$ is the convex cone in $\N(X)$ generated by
numerical classes of effective curves; $\overline{\NE}(X)$ is its
closure in $\N(X)$.

If $R$ is a half-line in $\N(X)$ and $D$ a $\Q$-Cartier divisor in
$X$, we will say that $D\cdot R>0$, $D\cdot R=0$, or $D\cdot R<0$, if
for any non zero element $\gamma\in R$ we have respectively
 $D\cdot \gamma>0$, $D\cdot \gamma=0$, or $D\cdot \gamma<0$.

If $\ph\colon X\to Y$ is a contraction, then the push-forward of
$1$-cycles gives a surjective linear map 
$$\ph_*\colon
\N(X)\la\N(Y),$$ 
and we set  $\NE(\ph):=\NE(X)\cap\ker\ph_*$.
We say that $\ph$ is \emph{elementary} if $\rho_X-\rho_Y=1$. 

Suppose that 
 $X$ is $\Q$-factorial and $\ph$ is elementary with $\dim\Exc(\ph)=n-1$.
Then $\Exc(\ph)$ is an irreducible divisor
and $\Exc(\ph)\cdot \NE(\ph)<0$.

For any irreducible
closed subset $Z$ of $X$, let $i\colon Z\hookrightarrow X$ be
the inclusion, and consider the push-forward of $1$-cycles
$i_*\colon\N(Z)\to\N(X)$. We define
$$\N(Z,X):=i_*(\N(Z))\subseteq\N(X).$$
Equivalently, $\N(Z,X)$ is the linear subspace of $\N(X)$ spanned by
classes of curves contained in $Z$.
Working with $\N(Z,X)$ instead of $\N(Z)$ means that we consider
curves in $Z$ modulo numerical equivalence in $X$, instead of
numerical equivalence in $Z$. Notice that $\dim\N(Z,X)\leq\rho_Z$.

\medskip

\noindent{\bf One-dimensional fibers in Mori contractions.}
The following Theorem collects results 
due to several people, see \cite[Lemma~2.12 and Th.~4.1]{AWaview} and
references 
  therein. Notice that  $X_0$
does not need to be complete.
\begin{thm}
\label{AW1}
Let $X_0$ be a smooth variety, $\ph_0\colon
X_0\to Y_0$ a Mori contraction, and $F$ a fiber of $\ph_0$ having a 
one-dimensional irreducible component $F_0$. Then $Y_0$ is smooth in
$\ph_0(F)$ and
either
$F=F_0\cong\pr^1$,
or $\ph_0$ is of fiber type and $F$ has two irreducible components,
both isomorphic to $\pr^1$.
\end{thm}
Suppose in particular that
 every fiber of $\ph_0$ has dimension at most $1$, so that
 $Y_0$ is smooth. If $\ph_0$ is of fiber type, 
we will say that $\ph_0$ is a \emph{conic bundle}. If
$\ph_0$ is birational, then it is the blow-up of a
smooth, codimension $2$ subvariety of $Y_0$; we will say that $\ph_0$
is of type $(n-1,n-2)^{sm}$.

Concerning the singular case, we have the following.
\begin{thm}[\cite{ishii}, Lemma 1.1]\label{ish}
Let $X$ be a projective variety with terminal singularities, and
$\ph\colon X\to Y$ an elementary birational Mori contraction with
fibers of dimension at most $1$. If $F_0$ is an irreducible component
of a non-trivial fiber of $\ph$, and $F_0$ contains a Gorenstein point
of $X$, then $F_0\cong\pr^1$ and $-K_X\cdot F_0\leq 1$.
\end{thm}

\noindent{\bf Fano varieties and Mori dream spaces.}
The notion of Mori dream space has been
introduced in \cite{hukeel}, where it is shown
that Fano $3$-folds are Mori dream spaces \cite[Cor.~2.16]{hukeel}. 
Moreover the authors 
 conjecture the same to hold in arbitrary
dimension. This has been confirmed recently in \cite{BCHM}, 
as an application of fundamental results on the minimal
model program.
\begin{thm}[\cite{BCHM}, Cor.~1.3.1]\label{tata}
Any smooth Fano variety is a Mori dream space.
\end{thm}
\noindent (In fact one can also allow singularities; here we
 consider only the smooth case.)

Being a Mori dream space
implies many important features with respect to Mori theory.
In the following remarks
we recall some consequences of Th.~\ref{tata} 
which will be used in the sequel.
\begin{remark}\label{asilo}
Let $X$ be a smooth Fano variety, and $X\dasharrow Y$ 
a ``rational contraction'' in the sense of
\cite{hukeel}. This means that there exists a normal and
$\Q$-factorial projective variety $X'$, and a factorization
$$X\dasharrow X'\la Y,$$
such that $X\dasharrow X'$ is an isomorphism in codimension $1$, and
$X'\to Y$ is a contraction.

Many well-known properties of $X$ hold for $Y$ too. The Mori cone
$\NE(Y)$ is closed and polyhedral.
For any contraction $\psi\colon Y\to Z$, $\NE(\psi)$ is a face of
$\NE(Y)$, which determines $\psi$ uniquely. 
 Conversely, for every face $F$ of $\NE(Y)$ there exists a
contraction $\psi$ of $Y$ such that $F=\NE(\psi)$.
Finally $\psi$ is elementary if and only if $\NE(\psi)$ has 
dimension one; we will call \emph{extremal ray} a $1$-dimensional face
of $\NE(Y)$.

This follows from the very definition of Mori dream space. Indeed
$X'$ is a ``small $\Q$-factorial
modification of $X$'', thus by \cite[Def.~1.10 and 
Prop.~1.11(2)]{hukeel} the properties above hold for $X'$. Then it is not
difficult to deduce the same for $Y$. 

\medskip

If $R=\NE(\psi)$ is an extremal ray of $\NE(Y)$, we say that $R$
is birational, divisorial, small, of fiber type, or of type $(a,b)$,
if the contraction $\psi$ is. Moreover we set $\Lo(R):=\Exc(\psi)$.

\medskip

Consider the special case where $\ph\colon X\to Y$ is an elementary
contraction. Then the extremal rays of $\NE(Y)$ are in bijection (via
$\ph_*$) with the $2$-dimensional faces of $\NE(X)$ containing the ray
$\NE(\ph)$, see \cite[2.5]{fanos}.
\end{remark}
\begin{remark}
\label{basic}
Let $Y$ be as in Rem.~\ref{asilo}, suppose moreover that it is
 $\Q$-factorial,
and consider a prime divisor $D\subset Y$.

There exists at least one extremal ray of
$\NE(Y)$ having positive intersection with $D$. Looking at the
associated contraction, one finds an elementary contraction
$$\psi\colon Y\la Z$$
such that $D\cdot\NE(\psi)>0$, in particular $D$ intersects every non
trivial fiber of $\psi$.

If $\psi$ is of fiber type, then $\psi(D)=Z$, hence
$$\psi_*\bigl(\N(D,Y)\bigr)=\N(Z)$$
and $\rho_Z\leq\dim\N(D,Y)$, $\,\rho_Y\leq\dim\N(D,Y)+1$.

If $\psi$ is birational, then $\Exc(\psi)\cap D\neq \emptyset$,
however $\Exc(\psi)\neq D$ (otherwise it should be
$D\cdot\NE(\psi)<0$), thus $\psi(D)\subset Z$ is a divisor.
We have two possibilities: either $\NE(\psi)\subset \N(D,Y)$ and
$\dim\N(\psi(D), Z)=\dim\N(D,Y)-1$, or $\NE(\psi)\not\subset
\N(D,Y)$ and
$\dim\N(\psi(D), Z)=\dim\N(D,Y)$. 
In this last case $\psi$ must be finite on $D$, hence every
non trivial fiber of $\psi$ is a curve.
\end{remark}
\begin{remark}\label{pera}
Let $X$ be a smooth Fano variety, and $D$ a prime 
divisor in $X$. By \cite[Prop.~1.11(1)]{hukeel}
there exists a finite sequence
\stepcounter{thm}
\begin{equation}\label{Dmmp}
X=X_0\stackrel{\sigma_0}{\dasharrow} X_1 \dasharrow\quad\cdots
\quad
\dasharrow X_{k-1}\stackrel{\sigma_{k-1}}{\dasharrow} X_k
\end{equation}
 where:
\begin{enumerate}[$\bullet$]
\item every $X_i$ is projective, normal, and $\Q$-factorial;
\item if $D_i\subset X_i$ is the proper transform of $D$, 
for $i=0,\dotsc,m-1$ there exists a birational
 extremal  ray $R_i$ of $\NE(X_i)$ such that
 $D_i\cdot R_i>0$, and
$\sigma_i$ is either the
  contraction of $R_i$ (if divisorial) or its flip (if small);
\item there exists an extremal ray of
  fiber type $R_k$
  of $\NE(X_k)$
  with $D_k\cdot R_k>0$. 
\end{enumerate}
See \cite[Def.~3.33 and~6.5]{kollarmori} for the definition of
flip. In the terminology of \cite{kollarmori,hukeel}
we are considering $(-D)$-flips, and
\eqref{Dmmp} is a Mori program for $-D$: since $D$ is effective, $-D$
can never become nef, so the program necessarily ends with a fiber
type contraction.
Notice that the choice of the extremal rays $R_i$'s is arbitrary among
the ones with positive intersection with $D_i$.
\end{remark}
\section{Divisors with Picard number $2$}\label{D}
Let $X$ be a smooth Fano variety and $D\subset
X$ a
 prime divisor. We recall that $\N(D,X)$ is the linear subspace of
$\N(X)$ spanned by classes of
curves contained in $D$, so that $\dim\N(D,X)\leq\rho_D$. 
The following result is proven in \cite{toru} under the assumption
that $\rho_D=1$, however the same proof works when  $\dim\N(D,X)=1$,
see  \cite[Prop.\
3.16]{fanos}.
\begin{proposition}[\cite{toru}, Prop.\ 5]\label{1}
Let $X$ be a smooth Fano variety of dimension $n\geq 3$, and $D\subset
X$ a
 prime divisor with $\dim\N(D,X)=1$. Then $\rho_X\leq 3$.
\end{proposition}
\noindent In particular we get $\rho_X\leq 3$
when $X$ has an elementary contraction of type
$(n-1,0)$.

In this section we consider the case where
 $\dim\N(D,X)=2$. Our
goal is to prove the following two results, which give a bound on
$\rho_X$ in some cases.
\begin{thm}\label{nuovo}
Let $X$ be a smooth Fano variety of dimension $n\geq 3$, and $D\subset
X$ a
 prime divisor with $\dim\N(D,X)=2$. Let $\ph\colon X\to Y$ be an
elementary 
contraction of $X$ with $D\cdot \NE(\ph)>0$. Then one of the following holds:
\begin{enumerate}[$(i)$]
\item $\rho_X=2$;
\item $\rho_X=3$ and $\ph$ is either
a conic bundle, or 
 of type $(n-1,0)$, or $(n-1,n-2)^{sm}$, or small;
\item $\ph$ is of type $(n-1,n-2)^{sm}$ and
  $\NE(\ph)\not\subset\N(D,X)$;
\item $\ph$ is
small, and there exists a smooth prime divisor
  $D'\subset X$, disjoint from $\Exc(\ph)$,
with a $\pr^1$-bundle structure, such that
 for any fiber $f$ we have 
 $D'\cdot f=-1$, $D\cdot f>0$, and $f\not\subset D$.
\end{enumerate}
\end{thm}
\noindent In the last case, we do not know whether the numerical class $[f]$
lies on an extremal ray of $\NE(X)$. However
$X$ is
the blow-up of a 
(possibly non projective)
complex manifold, in a smooth codimension $2$
subvariety, with exceptional divisor $D'$. 
\begin{lemma}\label{10/2}
Let $X$ be a smooth Fano variety of dimension $n\geq 3$, and
$D\subset X$ a prime divisor with $\dim\N(D,X)=2$. 

Suppose that
there exists an elementary divisorial contraction
$\ph\colon X\to Y$ such that
 $D\cdot\NE(\ph)=0$ and $\Exc(\ph)\cap D\neq\emptyset$.

Then either $\rho_X\leq 4$, or
there exists an extremal ray $R\neq\NE(\ph)$, of type $(n-1,n-2)^{sm}$, 
  such that $R\cdot \Exc(\ph)<0$ and $R+\NE(\ph)$ is a face of $\NE(X)$.
\end{lemma}
Notice that if $X$ is a toric Fano variety and $D\subset X$ is a prime
divisor which is closed with respect to the torus action, then
$\rho_X\leq 3+\dim\N(D,X)=3+\rho_D$ by \cite[Th.\ 2.4]{fano}; in
particular $\rho_X\leq 5$ when $\dim\N(D,X)=2$. However in general one
can not expect a similar bound, as the following example shows.
\begin{example}
Consider a Del Pezzo surface $S$ with $\rho_S=9$,
and let $X=S\times\pr^{n-2}$. Then $\rho_X=10$, and $X$ contains divisors 
$D= C\times\pr^{n-2}$, where $C\subset S$ is an irreducible curve, with
$\dim\N(D,X)=2$. 
\end{example}

\medskip

Before
proving Th.~\ref{nuovo} and Lemma~\ref{10/2}, we need some
preliminary properties. 
We fix a smooth Fano variety $X$ of dimension
$n\geq 3$ and a prime
divisor $D\subset X$, and
we carry out Mori's program for $-D$ as explained
in Rem.~\ref{pera}.  
We stop at $X_m$ when we get either
a contraction of fiber
type, or a birational extremal ray $R_m$ which is not contained in
$\N(D_m,X_m)$. 
Thus we obtain a sequence as \eqref{Dmmp}:
\stepcounter{thm}
\begin{equation}\label{mmp}
X=X_0\stackrel{\sigma_0}{\dasharrow} X_1 \dasharrow\quad\cdots
\quad
\dasharrow X_{m-1}\stackrel{\sigma_{m-1}}{\dasharrow} X_m
\end{equation}
 where moreover $R_i\subset\N(D_i,X_i)$ for
 $i=0,\dotsc,m-1$, and  
there exists an extremal ray $R_m$ of $\NE(X_m)$ with
$D_m\cdot R_m>0$, which is either of fiber type, or birational with
$R_m\not\subset\N(D_m,X_m)$.
\begin{lemma}
\label{26/6}
For every $i=0,\dotsc,m-1$ we have:
$$\dim\N(D_{i+1},X_{i+1})=\begin{cases}
\dim\N(D_i,X_i)-1 \ &\text{ if $R_i$ is divisorial;}\\
\dim\N(D_i,X_i) &\text{ if $R_i$ is small.}
\end{cases}$$
\end{lemma}
\begin{proof}
By construction we have $R_i\subset\N(D_i,X_i)$, thus
the statement is clear if $R_i$ is divisorial. 
Suppose that
 $R_i$ is small,  
let $\ph_i\colon
X_i\to Y_i$ be its contraction, $\ph'_{i}\colon X_{i+1}\to Y_i$
the flip of $\ph_i$,
and $R'_{i}:=\NE(\ph'_{i})$.
$$\xymatrix{
{X_{i}}\ar@{-->}[rr]^{\sigma_{i}}\ar[dr]_{\ph_i} &
&{X_{i+1}}\ar[dl]^{\ph_i'}\\
& {Y_i}&
}$$
Then $D_{i+1}\cdot
R'_{i}<0$ (see \cite[Cor.~6.4(4)]{kollarmori}), 
hence $R'_{i}\subset\N(D_{i+1},X_{i+1})$. This implies
the statement, because 
 $\ph_i(D_i)=\ph_i'(D_{i+1})$ and
$$\dim\N(D_i,X_i)=\dim\N(\ph_i(D_i),Y_i)+1=\dim\N(D_{i+1},X_{i+1}).$$
\end{proof}
\begin{corollary}\label{fibertype}
Suppose that in \eqref{mmp} the ray $R_m$ is of fiber type. Then
$$\rho_X\leq 1+\dim\N(D,X).$$
\end{corollary}
\begin{proof}
We have 
$$\rho_{X_{i+1}}=\begin{cases}
\rho_{X_i}-1 \ &\text{ if $R_i$ is divisorial;}\\
\rho_{X_i} &\text{ if $R_i$ is small.}
\end{cases}$$
Thus Lemma \ref{26/6} says that $\rho_{X_i}-\dim\N(D_i,X_i)$ is
constant, in particular
$\rho_{X}-\dim\N(D,X)=\rho_{X_m}-\dim\N(D_m,X_m)$.
If $R_m$ is of fiber type, then $\rho_{X_m}\leq 1+\dim\N(D_m,X_m)$
(see Rem.~\ref{basic}),
which gives the statement.
\end{proof}
Let $A_1\subset X_1$ be the
indeterminacy locus of $\sigma_0^{-1}$, and for $i\in\{2,\dotsc,m\}$
let
$A_i\subset X_i$ be the
union of the proper transform of $A_{i-1}\subset X_{i-1}$, with the
indeterminacy locus of $\sigma_{i-1}^{-1}$.
Then $X_i\smallsetminus A_i$ is isomorphic to an open subset of $X$,
and $$\Sing(X_i)\subseteq A_i\subset D_i.$$
Notice moreover that $\dim A_i>0$ whenever $R_{i-1}$ is small.
\begin{lemma}
\label{dentini}
Let $i\in\{1,\dotsc,m\}$, and assume that $-K_{X_j}\cdot R_j>0$
for every $j=0,\dotsc,i-1$. 
Then $X_1,\dotsc,X_i$ have terminal singularities.
Moreover if $C\subset X_i$ is an irreducible curve not contained in
$A_i$, and $C_0\subset X$ its proper transform, we have
$$-K_{X_i}\cdot C\geq -K_X\cdot C_0,$$
with strict inequality whenever $C\cap A_i\neq \emptyset$.
\end{lemma}
\begin{proof}
We assume that the statement holds for $i-1$, and
consider $\sigma_{i-1}\colon X_{i-1}\dasharrow X_{i}$.
Suppose that $\sigma_{i-1}$ is a flip, and consider a common
resolution of $X_{i-1}$ and $X_i$:
$$\xymatrix{
&{\widehat{X}}\ar[ld]_f\ar[rd]^g & \\
{X_{i-1}}\ar@{-->}[rr]_{\sigma_{i-1}} & &{X_{i}}
}$$
Let $G_1,\dotsc,G_r\subset \widehat{X}$ be the exceptional
divisors, and write
$$K_{\widehat{X}}=f^*(K_{X_{i-1}})+\sum_{j=1}^ra_jG_j=g^*(K_{X_{i}})
+\sum_{j=1}^rb_jG_j,\quad\text{with }a_j,b_j\in\Q.$$
Since $X_{i-1}$ has terminal singularities and
$-K_{X_{i-1}}\cdot R_{i-1}>0$, we have $b_j\geq a_j>0$ for every
$j=1,\dotsc,r$ by \cite[Lemma
  3.38]{kollarmori}, thus also
$X_i$ has terminal singularities. 

The curve $C\subset X_i$ is not contained in $A_i$, 
hence $C$ intersects the
open subset where $X_{i-1}$ and $X_i$ are non singular and isomorphic.
If $\widetilde{C}\subset X_{i-1}$
and $\widehat{C}\subset \widehat{X}$ are the proper transforms of $C$,
then  $G_j\cdot \widehat{C}\geq 0$ for every $j$, and
we get
$$-K_{X_{i}}\cdot C=-K_{X_{i-1}}\cdot \widetilde{C}
+\sum_{j=1}^r(b_j-a_j)G_j\cdot \widehat{C}\geq -K_{X_{i-1}}\cdot
\widetilde{C} 
\geq -K_X\cdot C_0.$$
Now suppose that $C\cap A_i\neq\emptyset$. If $\widetilde{C}\cap
A_{i-1}\neq\emptyset$, then $-K_{X_{i-1}}\cdot \widetilde{C}> -K_X\cdot C_0$. 
Otherwise $\widetilde{C}$ must intersect $\Lo(R_{i-1})$, thus
there exists some $j_0$ such that $f(G_{j_0})\subseteq\Lo(R_{i-1})$ and
$G_{j_0}\cdot \widehat{C}>0$. Again by \cite[Lemma
  3.38]{kollarmori} we have $b_{j_0}>a_{j_0}$, therefore
$-K_{X_{i}}\cdot C>-K_{X_{i-1}}\cdot \widetilde{C}$ and we are done.

The case where $\sigma_{i}$ is a divisorial contraction is similar and
shorter.
\end{proof}
\begin{lemma}
\label{23/6}
Let $i\in\{1,\dotsc,m\}$, and assume that $-K_{X_j}\cdot R_j>0$
for every $j=0,\dotsc,i-1$. 

Consider a
birational elementary contraction
$\psi\colon X_i\to Y$ be  such that $D_i\cdot\NE(\psi)>0$ and
$\NE(\psi)\not\subset\N(D_i,X_i)$. 

Then $\Exc(\psi)$ is disjoint
from $A_i$, and 
$\psi_{|X_i\smallsetminus A_i}$ is a Mori contraction of type $(n-1,n-2)^{sm}$.
\end{lemma}
\begin{proof}
Let $F$ be a non trivial fiber of $\psi$. Then
$F$ must meet $D_i$, 
on the other hand $\psi$ is finite on $D_i$. Thus 
$F$ is a curve which intersects $D_i$ in finitely many
points, in particular $F$ can not be contained in $A_i$.

Using Lemma \ref{dentini}
 we see that $-K_{X_i}\cdot F>0$, namely
$\psi$ is a Mori contraction; moreover
$\dim(F\cap\Sing(X_i))\leq 0$. We can now apply 
Th.~\ref{ish} to deduce that  $-K_{X_i}\cdot F_0\leq 1$
for any irreducible component $F_0$ of
 $F$.

Again by Lemma \ref{dentini}, this shows that $F$ can not intersect $A_i$; 
in particular $\Exc(\psi)$ is contained in the
smooth locus of $X_i$, and the statement follows.
\end{proof}
\begin{lemma}
\label{24/6}
Let $i\in\{1,\dotsc,m\}$, and assume that $-K_{X_j}\cdot R_j>0$
for every $j=0,\dotsc,i-1$. 

If $\dim\N(D_i,X_i)=1$
and $\dim A_i>0$, 
then $i=m$,
$\rho_{X_m}\leq 2$, and every elementary contraction $\psi\colon
X_m\to Y$ such that  $D_m\cdot \NE(\psi)>0$ is of fiber type.
In particular $R_m$ is of fiber type.
\end{lemma}
\begin{proof}
Let $\psi\colon X_i\to Y$
 be an elementary contraction such that $D_i\cdot \NE(\psi)>0$.

If  $\psi$ is birational, 
it can not be $\NE(\psi)\subset\N(D_i,X_i)$, otherwise $\psi(D_i)$ is
a point and $D_i=\Exc(\psi)$, which contradicts $D_i\cdot
\NE(\psi)>0$.
 On the other hand,
if $\NE(\psi)\not\subset\N(D_i,X_i)$, then 
 Lemma \ref{23/6} implies
that $\Exc(\psi)$ is a divisor disjoint from $A_i$. 
However this is again impossible, because
$\Exc(\psi)\cap D_i\neq\emptyset$, so there exists some curve
$C\subset D_i$ with $\Exc(\psi)\cdot C>0$. Since all curves in $D_i$
are numerically proportional, the same must hold for every curve 
$C\subset D_i$. Now choosing $C\subseteq A_i$ we get $\Exc(\psi)\cap
A_i\neq\emptyset$.

Therefore $\psi$ is of fiber type, $\rho_{X_i}\leq 2$, and it must be
$i=m$. 
\end{proof}
\begin{proof}[Proof of Th.~\ref{nuovo}]
We assume that $\rho_X\geq 3$, and show that one of $(ii)$, $(iii)$,
or $(iv)$ holds. 
Let's consider the possibilities for $\ph$.

If $\ph$ is of fiber type, then $\ph(D)=Y$, so $\rho_X=3$ and
$\ph$ is finite on $D$. Then $\ph$ must have only $1$-dimensional
fibers, it is a conic bundle, and we are in $(ii)$.

\medskip

Suppose now that $\ph$ is birational. 

If $\ph$ is of type $(n-1,0)$, then $\rho_X=3$ by
Prop.~\ref{1}, so we are again in $(ii)$.
If instead
$\NE(\ph)\not\subset\N(D,X)$, then $\ph$ is of type $(n-1,n-2)^{sm}$
and we are in $(iii)$.

Thus we assume that $\ph$ is not of type $(n-1,0)$ and that
$\NE(\ph)\subset\N(D,X)$.

\medskip

Consider the sequence \eqref{mmp}. We can assume that $R_0=\NE(\ph)$, so
that $m\geq 1$. Then $\dim A_1>0$, because if $\ph$ is divisorial then
$A_1=\ph(\Exc(\ph))$. 

\medskip

Suppose that $R_m$ is of fiber type. Then Cor.~\ref{fibertype}
gives $\rho_X=3$, and in order to get $(ii)$ we are left to show that
$\ph$ is either small or of type $(n-1,n-2)^{sm}$.

Let's assume that $\ph$ is divisorial.
Then $\dim\N(D_1,X_1)=1$, and
Lemma~\ref{24/6} yields that $m=1$, namely
we have:
$$X\stackrel{\ph}{\la} X_1\stackrel{\psi}{\la}Y,$$
where $\psi$ is the contraction of $R_1$, and is of fiber type.

We have
$\rho_Y=1$, so $Y$ is not a point.
Since 
all curves contained in $D_1$ are numerically proportional, $\psi$
must be finite on $D_1$. Then
  every fiber of $\psi$ has dimension $1$ and $\dim Y=n-1$.

Notice that $\psi$ is finite
on $A_1=\ph(\Exc(\ph))$, because $A_1\subset
D_1$. Choose a point $x_1\in A_1$. The fiber
$\psi^{-1}(\psi(x_1))$ has dimension $1$ and is not contained in
$A_1$, hence $\ph^{-1}(\psi^{-1}(\psi(x_1))$ has some
$1$-dimensional irreducible component. Then Th.~\ref{AW1} applied to
$\psi\circ\ph$ yields that $\ph^{-1}(\psi^{-1}(\psi(x_1))$ is
$1$-dimensional and has exactly two irreducible components. This means
that $\psi^{-1}(\psi(x_1))\cap A_1=\{x_1\}$
(i.e.\ $\psi$ is injective on $A_1$), and the two
components are $\ph^{-1}(x_1)$
and the proper transform of $\psi^{-1}(\psi(x_1))$. 

Therefore every non trivial fiber of $\ph$ is $1$-dimensional, so
$X_1$ is smooth and $\ph$ is of type $(n-1,n-2)^{sm}$. 
In fact it is not difficult to show that $X_1$ is Fano and
that $\psi$
is a smooth morphism.

\medskip

Let's consider now the case where $R_m$ is birational, and 
show that this gives~$(iv)$.

We claim that $-K_{X_i}\cdot
R_i>0$ for every $i=0,\dotsc,m$. Indeed this is true for $i=0$. 
Fix $i\in\{1,\dotsc,m\}$ and
assume
that $-K_{X_j}\cdot
R_j>0$ for $j=0,\dotsc,i-1$. 

We observe that $\dim A_i>0$. 
This is clear if $i=1$ or
if $R_{i-1}$ is small. Suppose that 
$i>1$ and that $R_{i-1}$ is divisorial, so that $\sigma_{i-1}$ is its
contraction. Since $\dim\N(D,X)=2$, it follows from Lemma \ref{26/6}
that there is at most one divisorial ray among $R_0,\dotsc,R_{m-1}$.
Thus $R_{i-2}$ is small, and $A_{i-1}$
contains the indeterminacy locus $L$ of $\sigma_{i-2}^{-1}$, which is the
locus of a small extremal ray of $\NE(X_{i-1})$. Then $\sigma_{i-1}$
is finite on $L$ and $\sigma_{i-1}(L)\subset A_i$ has positive
dimension.

Therefore  Lemma \ref{24/6} implies that
$\dim\N(D_i,X_i)>1$, thus $\dim\N(D_i,X_i)=2$, and
$R_{i-1}$ is small.

Let $R'_{i-1}$ be the small 
extremal ray of $\NE(X_i)$ whose contraction is the flip of $R_{i-1}$ 
in $X_{i-1}$. Then
 $-K_{X_i}\cdot R'_{i-1}<0$ and $D_i\cdot R'_{i-1}<0$
(see \cite[Cor.~6.4(4)]{kollarmori}), thus
$R'_{i-1}\subset\N(D_i,X_i)$ and
$$\N(D_i,X_i)\cap\NE(X_i)=R_i+R'_{i-1}.$$
Since by Lemma \ref{dentini} 
the divisor $D_i$ contains curves of positive anticanonical degree, and
$-K_{X_i}\cdot R'_{i-1}<0$, we must have $-K_{X_i}\cdot R_i>0$.
We have also shown
that $R_0,\dotsc,R_{m-1}$ are small, in particular $\ph$ is small.

Now it follows from Lemma \ref{23/6} that
 $\Lo(R_m)\cap A_m=\emptyset$, 
and that $R_m$ is of type
$(n-1,n-2)^{sm}$. Therefore the proper transform of $\Lo(R_m)$ in $X$
yields a divisor $D'$ as in $(iv)$, and we are done.
\end{proof}
We need one more Lemma before proving Lemma \ref{10/2}.
\begin{lemma}\label{3/4}
Let $X$ be a smooth Fano variety of dimension $n$, and $\ph_1\colon
X\to Y_1$ a divisorial elementary contraction. Let $\psi\colon Y_1\to
Z$ be an elementary birational contraction with fibers of dimension at
most $1$.

Consider the elementary contraction
 $\ph_2\colon X\to Y_2$ such that
 $\NE(\psi\circ\ph_1)=\NE(\ph_1)+\NE(\ph_2)$, and set
$E_i:=\Exc(\ph_i)\subset X$ for $i=1,2$.
$$\xymatrix{
X\ar[d]_{\ph_1}\ar[r]^{\ph_2}& {Y_2}\ar[d] \\
{Y_1}\ar[r]_{\psi} & Z}$$
Then $Y_2$ is smooth, $\ph_2$ is of type $(n-1,n-2)^{sm}$,
 and $\Exc(\psi)=\ph_1(E_2)$. Moreover one of the
following holds:
\begin{enumerate}[$(i)$]
\item $\psi$ is a divisorial Mori contraction, 
$\Exc(\psi)\cap
  \ph_1(E_1)$ is a union of fibers of $\psi$, $E_1\cdot\NE(\ph_2)=0$,
  and $E_1\neq E_2$;
\item $\psi$ is small,   $\Exc(\psi)=\ph_1(E_1)$,
  $E_1\cdot\NE(\ph_2)<0$, and $E_1=E_2$.
\end{enumerate}
\end{lemma}
\begin{proof}
Let $F$ be a non trivial fiber of $\psi$, then
$(\ph_1)^{-1}(F)$ is a fiber of $\psi\circ\ph_1\colon X\to
Z$. By Th.~\ref{AW1}, if $(\ph_1)^{-1}(F)$ has an irreducible
component of dimension $1$, then $(\ph_1)^{-1}(F)\cong\pr^1$. 
This means that either
$F\subseteq\ph_1(E_1)$ or $F\cap\ph_1(E_1)=\emptyset$.
Therefore $\Exc(\psi)\cap\ph_1(E_1)$ is a union of fibers
of $\psi$.

\medskip

Now let $F'$ be a non trivial fiber of $\ph_2$. Then $\ph_1(F')$ is
contained in a non trivial fiber of $\psi$, thus
$\ph_1(F')\subseteq\Exc(\psi)$
and $\dim\ph_1(F')=1$. But $\ph_1$ is finite on $F'$, so $\dim F'=1$,
and $\ph_2$ is birational with fibers of dimension at most $1$. Thus
 $Y_2$ is smooth and $\ph_2$ is of type $(n-1,n-2)^{sm}$.
We also have $\ph_1(E_2)\subseteq\Exc(\psi)$.

Notice  that if $F'$ intersects $E_1$, then $\ph_1(F')$ intersects
$\ph_1(E_1)$, hence $\ph_1(F')\subseteq\ph_1(E_1)$. 

\medskip

Suppose that $E_1\neq E_2$. Then $\ph_1(E_2)$ is a divisor contained
in $\Exc(\psi)$,
hence $\psi$ is divisorial
 with exceptional locus
$\ph_1(E_2)$ (notice that $Y_1$ is $\Q$-factorial because $\ph_1$ is
divisorial).
Since $\ph_1(E_2)$
can not be contained in $\ph_1(E_1)$, $\psi$
has non trivial fibers which are disjoint from $\ph_1(E_1)$, so
it is
a Mori contraction. 
Moreover there must be fibers of $\ph_2$ which are disjoint from
$E_1$, hence $E_1\cdot\NE(\ph_2)=0$.

\medskip

Assume that $E_1=E_2$, so that $E_1\cdot \NE(\ph_2)<0$. Clearly the
exceptional locus of $\psi\circ\ph_1$ contains $E_1$. On the other
hand, every curve in $\NE(\psi\circ\ph_1)=\NE(\ph_1)+\NE(\ph_2)$ has
                                negative intersection with $E_1$,
                                hence it is contained in $E_1$, namely 
$\Exc(\psi\circ\ph_1)=E_1$. This yields $\Exc(\psi)=\ph_1(E_1)$, and
                                $\psi$ is small.
\end{proof}
\begin{proof}[Proof of Lemma \ref{10/2}]
We set $E:=\Exc(\ph)\subset X$, $A:=\ph(E)\subset Y$, and
$D_Y:=\ph(D)\subset Y$. Notice that there are non trivial fibers of
$\ph$ disjoint from $D$, and others contained in $D$. Therefore 
we have $\dim\N(D_Y,Y)=1$, but $A$ is not contained in $D_Y$: this is
the main difference with respect to the situation of Lemma \ref{24/6}.

Let $\psi\colon Y\to Z$
 be an elementary contraction of $Y$ with 
$D_Y\cdot \NE(\psi)>0$, as in Rem.~\ref{basic}.

If $\psi$ is of fiber type, we have $\rho_Z\leq 1$ and $\rho_X\leq 3$.
 
Suppose that $\psi$ is birational. As in the proof of Lemma \ref{24/6}
we see that $\psi$ is finite on
$D_Y$ and its fibers have dimension at most $1$, thus Lemma \ref{3/4}
applies; in particular $\Exc(\psi)\cap   A$ is a union of
fibers of $\psi$.

\medskip

If $\psi$ is not divisorial, then
Lemma \ref{3/4} $(ii)$ gives an extremal ray $R$ as in the statement. 
If $\psi$ is divisorial, we
are in Lemma \ref{3/4} $(i)$, thus $\psi$ is a Mori contraction and
$Z$ is $\Q$-factorial.

Set $D_Z:=\psi(D_Y)\subset Z$. Then $D_Z$ is a prime divisor in $Z$ with
$\dim\N(D_Z,Z)=1$, and $D_Z\supset\psi(\Exc(\psi))$. 

Let $\xi\colon Z\to W$  be an elementary contraction of $Z$ with
$D_Z\cdot\NE(\xi)>0$, as in Rem.~\ref{basic}.
If $\xi$ is of fiber type, we get $\rho_W\leq 1$ and
$\rho_X\leq 4$.

\medskip

Suppose that $\xi$ is birational; as before
it is finite over $D_Z$ and has
fibers of dimension at most $1$. 
Set $\eta:=\xi\circ\psi$, and let $\psi_1\colon Y\to Z_1$ be the
elementary contraction of $Y$
such that $\NE(\eta)=\NE(\psi)+\NE(\psi_1)$:
$$\xymatrix{
X\ar[r]^{\ph} &Y\ar[dr]^{\eta}\ar[d]_{\psi}\ar[r]^{\psi_1}& {Z_1}\ar[d] \\
&{Z}\ar[r]_{\xi} & W}$$
Again, $\psi_1$ is birational with fibers of dimension at most $1$,
thus Lemma \ref{3/4} applies. 
Either $\psi_1$ is not divisorial and we
get again an extremal ray $R$ as in the statement, 
or $\psi_1$ is a divisorial Mori contraction, and
$\Exc(\psi_1)\cap A$ is a union of fibers of $\psi_1$. We
show that this last case leads to a contradiction.

Every curve in $\NE(\eta)$ has positive
anticanonical degree, thus $\eta$ is a Mori contraction.

If $\Exc(\psi)=\Exc(\psi_1)$, then every curve in $\NE(\eta)$ has
negative intersection with $\Exc(\psi)$, hence $\Exc(\eta)=\Exc(\psi)$
and $\Exc(\xi)=\psi(\Exc(\psi))$. However this is impossible, because
$\xi$ is finite on $D_Z$ which contains $\psi(\Exc(\psi))$.

Therefore $\Exc(\psi)\neq\Exc(\psi_1)$. Then $\psi(\Exc(\psi_1))$ is a
divisor contained in $\Exc(\xi)$, which means that
$\Exc(\xi)=\psi(\Exc(\psi_1))$ and $\xi$ is divisorial. As 
in the proof of Lemma \ref{24/6}, we
see that $\Exc(\xi)$ must intersect every curve contained in $D_Z$,
and $\dim\psi(\Exc(\psi))=n-2\geq 1$,
hence $\Exc(\xi)\cap \psi(\Exc(\psi))\neq\emptyset$.
Then $\dim
(\psi(\Exc(\psi))\cap\Exc(\xi))\geq n-3$
and since $\xi$ is finite on $\psi(\Exc(\psi))$, we get
$$\dim\xi\bigl(\psi(\Exc(\psi))\cap\Exc(\xi)\bigr)\geq n-3.$$
We claim that
\stepcounter{thm}
\begin{equation}\label{sing}
\dim\xi\bigl(\psi(\Sing(Y))\cap\Exc(\xi)\bigr)\leq n-4.
\end{equation}
First let's see that \eqref{sing} allows to conclude the proof. 
Since both $\psi$ and $\xi$ are Mori contractions with fibers of
dimension at most $1$, we have $\Sing(W)\subseteq\eta(\Sing(Y))$. Thus
\eqref{sing} implies that there exists a point $w_0\in
W\smallsetminus\Sing(W)$ such that the fiber $\xi^{-1}(w_0)$ has
dimension $1$ and intersects $\psi(\Exc(\psi))$. 
Now restricting $\eta$ to a contraction
$Y\smallsetminus\eta^{-1}(\Sing(W))\to W\smallsetminus\Sing(W)$, 
 we can apply Th.~\ref{AW1}  to $\eta^{-1}(w_0)$
as in the proof of Lemma \ref{3/4} and get a
contradiction. 

\medskip

Let's show \eqref{sing}. If $\dim A=n-2$, then $\dim\Sing(Y)\leq n-4$,
so  \eqref{sing} holds.
If $\dim A\leq n-3$, we still have $A\supseteq\Sing(Y)$, 
thus it is enough to show that $\dim\xi(\psi(A)\cap\Exc(\xi))\leq n-4$.
This is clear if 
$\psi(A)$ is not contained in $\Exc(\xi)$.
If instead
$\psi(A)\subseteq\Exc(\xi)=\psi(\Exc(\psi_1))$, we get
$$A=\psi^{-1}(\psi( A))\subseteq
\psi^{-1}(\Exc(\xi))\subseteq
\Exc(\psi)\cup \Exc(\psi_1).$$
 Since $ A$ is irreducible, it is contained either in 
$\Exc(\psi)$, or in $\Exc(\psi_1)$, and it is a union of fibers of
both $\psi$ and $\psi_1$.
In any cases we get $\dim\eta( A)\leq n-4$, and we are done.
\end{proof}
\section{Elementary contractions of type $(n-1,1)$}
\label{31}
Throughout this section, we fix the following notation:

\smallskip

\stepcounter{thm}
\noindent {\bf(\thethm)} \label{setup}
$X$ is a
smooth Fano variety of dimension $n\geq 4$, and $R_1$ is an extremal
ray of type $(n-1,1)$. For any integer $i\in\mathbb{Z}_{\geq 0}$, if $R_i$
is an extremal ray  of $\NE(X)$, we denote by $\ph_i\colon X\to
Y_i$ the associated contraction, and we set $E_i:=\Exc(\ph_i)$.

\smallskip

Our goal is to bound $\rho_X$;
notice that $\rho_X\geq 2$ by our assumptions.

We observe first of all that
since 
$\ph_1(E_1)$ is a curve, we have
$\dim\N(\ph_1(E_1),Y_1)=1$ and
$$\dim\N(E_1,X)=2,$$
thus we can apply to $E_1$ the results of the preceding section. Indeed
there exists some extremal
ray $R_2$ with $E_1\cdot R_2>0$, and by Th.~\ref{nuovo}
we can conclude at once that
$\rho_X\leq 3$ unless $R_2$ is small, or of type $(n-1,n-2)^{sm}$.
More precisely, we show the following.
\begin{thm}\label{positive}
Let $X$ and $R_1$ be as in $(\ref{setup}.1)$, and
let $R_2$ be an extremal ray with $E_1\cdot R_2>0$. Then one of the
following holds.
\begin{enumerate}[$(i)$]
\item $\rho_X\leq 4$, more precisely we have the possibilities:

$\quad$$\ph_2$ is of type $(n,n-1)$, $(n,n-2)$, or
 $(n-1,n-3)$, and $\rho_X=2$;

$\quad$$\ph_2$ is a conic bundle and $\rho_X=3$;

$\quad$$\ph_2$ is of type $(n-2,n-4)$ and $\rho_X\leq 3$;

$\quad$$n=4$, $\ph_2$ is of type $(2,0)$, and $\rho_X=4$;

$\quad$$\ph_2$ is of type $(n-1,n-2)$ and $\rho_X\leq 4$.

\item $\ph_2$  is of type $(n-1,n-2)^{sm}$, $E_2\cdot R_1=0$, and
  there exists an extremal ray $R_0\neq R_1$ 
such  that $E_1\cdot R_0<0$.
\end{enumerate}
\end{thm}
\noindent Case $(ii)$ will be treated in Prop.~\ref{11/4}, 
where we will show that $\rho_X\leq 5$.
\begin{remark}[Classification results by T.\ Tsukioka]\label{tt}
Suppose that $Y_1$ is smooth and $\ph_1$ is the blow-up of a smooth curve.
When $\ph_2$ is of type $(n,n-2)$, 
the possible $X$ and $Y_1$ are classified in
\cite{toruCRAS}. 
Moreover if $n=4$, $\ph_2$ is of type $(3,1)$, and $E_2$ is smooth,
then it is shown in \cite{toru4fold} that $Y_1\cong\pr^4$ and
$\ph_1(E_1)$ is an elliptic curve of degree $4$ in $\pr^4$.
\end{remark}
\begin{proof}[Proof of Th.~\ref{positive}]
First notice that every non trivial
fiber $F$ of
$\ph_2$ has dimension at most $2$. In fact $F\cap E_1\neq\emptyset$
and $\ph_1$ is finite on it, so that
$$\dim F-1\leq \dim F\cap E_1=\dim \ph_1(F\cap E_1)\leq\dim\ph_1(E_1)=1.$$
This, together with Th.~\ref{nuovo}, implies the statement, unless
we are in cases $(iii)$ or $(iv)$ of Th.~\ref{nuovo}.

We consider first case $(iii)$, so
we assume that $\ph_2$ is of type $(n-1,n-2)^{sm}$
and 
$R_2\not\subset\N(E_1,X)$.  
We will distinguish the two cases $E_2\cdot R_1=0$
and $E_2\cdot R_1>0$.

\medskip

Suppose that $E_2\cdot R_1=0$. 
 Then $E_2$ must contain some fiber $F'$ of
$\ph_1$ of dimension $n-2$. Since $\ph_2$ is finite on $F'$, we have
$\ph_2(F')=\ph_2(E_2)$, hence 
$$\N(E_2,X)=\R R_2+\N(F',X)=\R(R_1+R_2).$$ 
Then Lemma
\ref{10/2} applies to $E_2$ and $\ph_1$, and yields that 
either $\rho_X\leq 4$, or
we have $(ii)$. See Rem.~\ref{vaccinazione} for a
more precise description of this case.

\medskip

\begin{center}
\scalebox{0.50}{\includegraphics{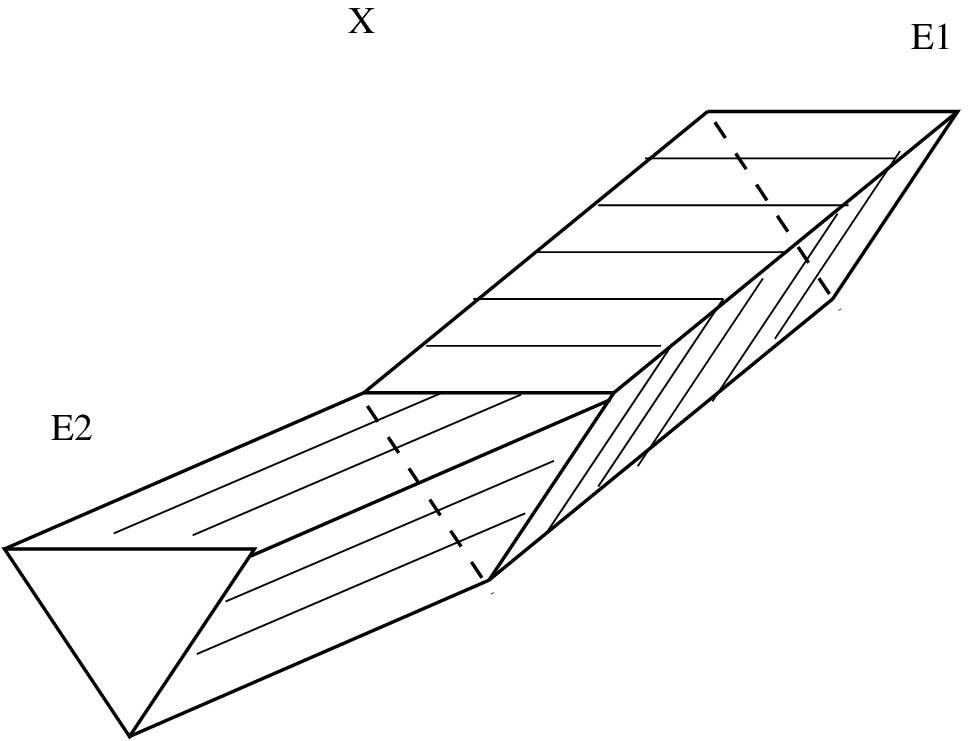}}

\medskip

{\footnotesize Figure 1: the case $R_2\not\subset\N(E_1,X)$ and
  $E_2\cdot R_1=0$.} 
\label{fig1}
\end{center}

\medskip

Assume now that $E_2\cdot R_1>0$, and consider $D:=\ph_2(E_1)\subset
Y_2$ and $A:=\ph_2(E_2)\subset D$. Then $A$ is smooth of dimension
$n-2$. We observe that if $C\subset Y_2$ is an irreducible curve not
contained in $A$, then $-K_{Y_2}\cdot C\geq 1$, with strict inequality
whenever $C\cap A\neq\emptyset$, by  Lemma \ref{dentini}.

\medskip 

\begin{center}
\scalebox{0.50}{\includegraphics{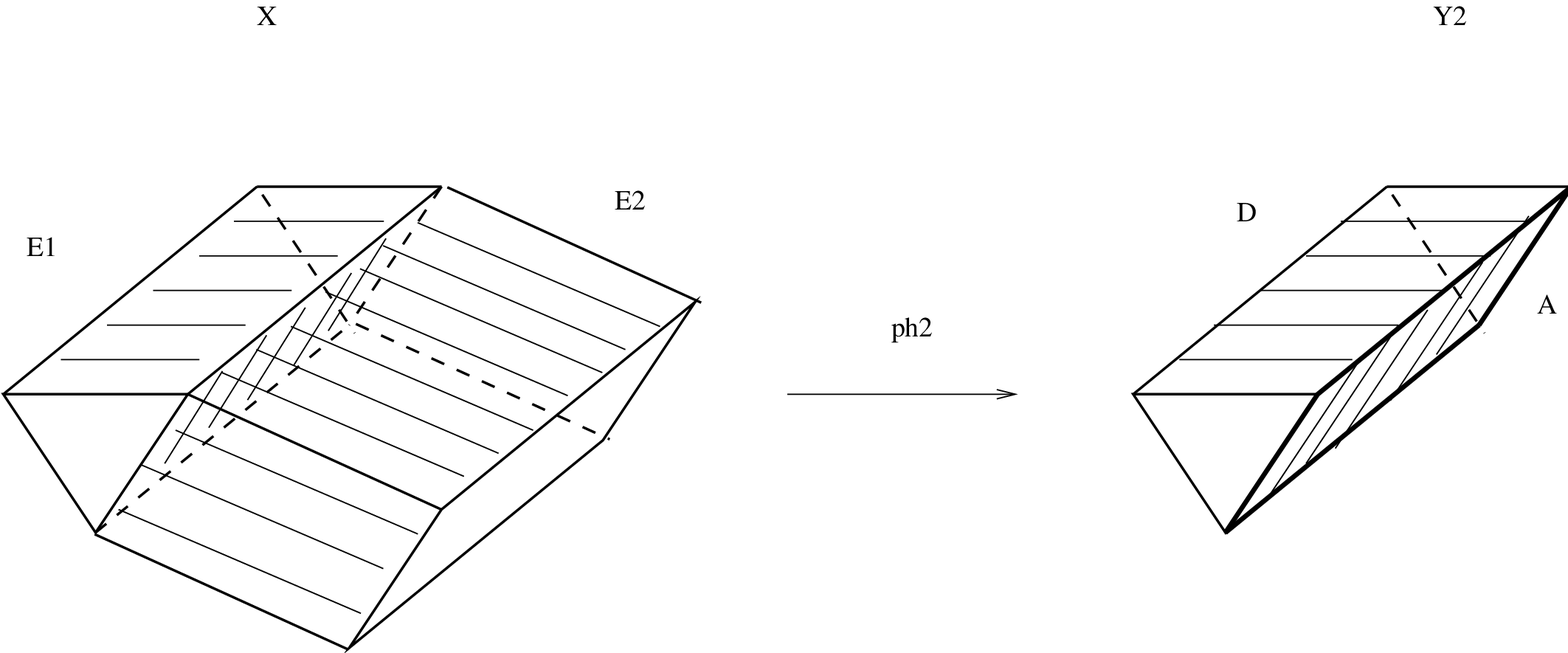}}

\medskip

{\footnotesize Figure 2: the case $R_2\not\subset\N(E_1,X)$ and
  $E_2\cdot R_1>0$.} 
\end{center}

\medskip

We first suppose that $Y_2$ is Fano, and 
 apply Th.~\ref{nuovo} to $D\subset Y_2$
(just choosing any extremal ray of $\NE(Y)$ having positive
intersection with $D$). If
$\rho_{Y_2}\leq 3$, then $\rho_X\leq 4$. 
Otherwise we are in cases $(iii)$ or $(iv)$ of 
Th.~\ref{nuovo}, and
there exists a smooth prime divisor $D'\subset Y_2$,
with a $\pr^1$-bundle structure, such that for any fiber $f$ we have
$-K_{Y_2}\cdot f=1$, $f\cdot D>0$, and $f\not\subset D$. 
Then $f\cap A=\emptyset$, namely
$D'$ can not intersect $A$, so its
inverse image
 $D''$ in $X$ is a prime divisor which intersects $E_1$ but
is disjoint from $E_2$. This is impossible, because either $D''\cdot
R_1=0$ and $D''$ contains some non trivial fiber of $\ph_1$, or
$D''\cdot R_1>0$ and $D''\cap C\neq\emptyset$ for some irreducible
curve $C\subset E_2$ with $[C]\in R_1$.

\medskip

Assume that $Y_2$ is not Fano. This means that there exists
some extremal ray of $\NE(Y_2)$ with non positive anticanonical
degree. Let's consider the associated contraction
$$\widetilde{\psi}\colon Y_2\la\widetilde{Z},$$
and notice that $\Exc(\widetilde{\psi})\subseteq
A\subset D$. Then any non trivial fiber of $\widetilde{\psi}$ must be
$1$-dimensional. In fact 
 if  $\widetilde{\psi}$
had a fiber $F$ with $\dim F\geq 2$, then we would have 
$$\dim(\ph_2^{-1}(F)\cap E_1)\geq 2,$$
thus $R_1\subset\N(\ph_2^{-1}(F),X)$ and
$(\ph_2)_*(R_1)\subset\N(F,Y_2)$. This implies that $\NE(\widetilde{\psi})=
(\ph_2)_*(R_1)$, which is impossible because $\Exc(\widetilde{\psi})$
should contain all $D$.

Therefore $\widetilde{\psi}$ is small with fibers of dimension at most
$1$. By Lemma \ref{3/4} we see that
there exists an extremal ray $R_3$
of $\NE(X)$ of type
$(n-1,n-2)^{sm}$ 
such that $E_2\cdot R_3<0$ and $R_2+R_3$ a face of $\NE(X)$; in
particular $\rho_X\geq 3$.

\medskip

\begin{center}
\scalebox{0.50}{\includegraphics{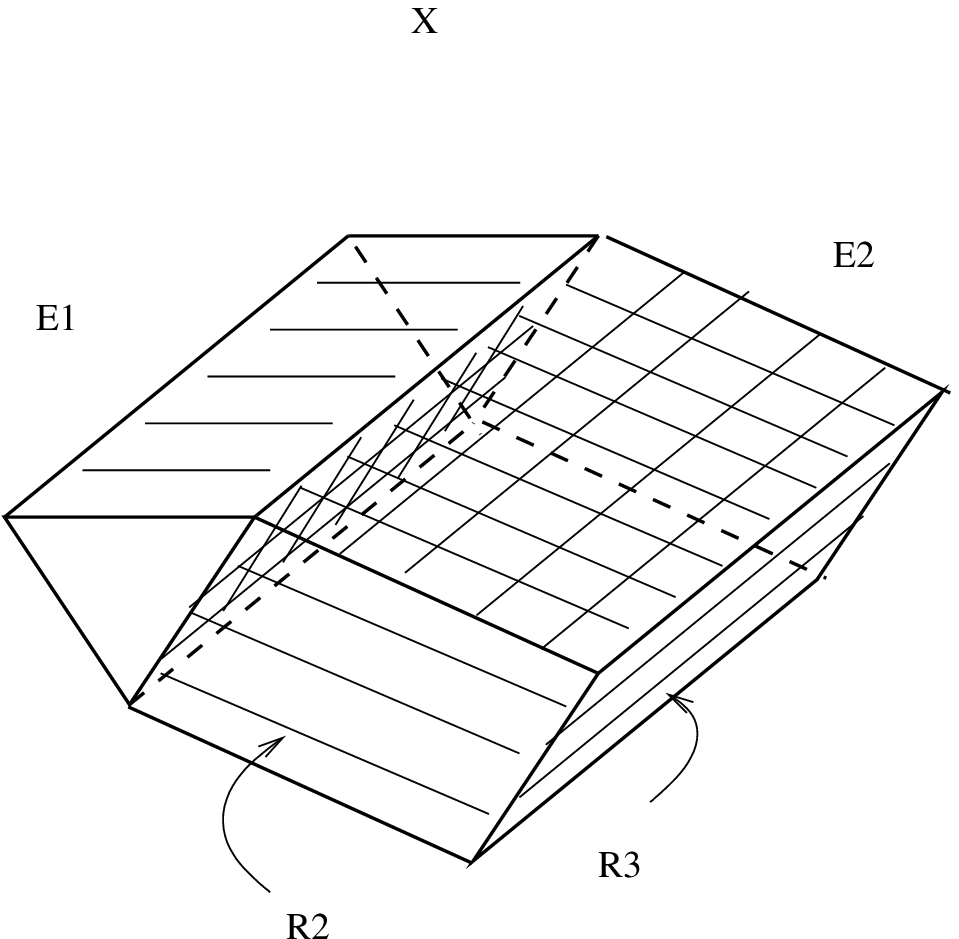}}

\medskip

{\footnotesize Figure 3: the case $Y_2$ not Fano.}

\end{center}

\medskip

We need to make some considerations on $E_2$, in order to show the
following:
\stepcounter{thm}
\begin{equation}
\label{vernavola}
\text{$E_1\cdot R_3=0$, and for every curve $C\subset E_2$ we have
$[C]\in R_1+R_2+R_3$.}
\end{equation}
Observe first of all
that $E_2$ is smooth, and $\ph_{2|E_2}$ and $\ph_{3|E_2}$ are
$\pr^1$-bundles. 
 Moreover we have $\N(E_2,X)=\R( R_1+R_2+R_3)$. Indeed $R_1$,
$R_2$ and $R_3$ are contained in $\N(E_2,X)$. On the other hand since
$E_1$ meets every fiber of $\ph_{2|E_2}$ we have
$\ph_2(E_1\cap E_2)=\ph_2(E_2)$, so that 
$$\R( R_1+R_2+R_3)
\subseteq\N(E_2,X)\subseteq \R R_2+\N(E_1,X).$$
However $\dim(\R R_2+\N(E_1,X))=3$, so the
inclusions above are equalities.

Let $T$ be the normalization of $\ph_1(E_2)$,
 $\xi\colon E_2\to T$ the contraction induced by
$(\ph_1)_{|E_2}$, and $i\colon E_2\hookrightarrow X$ the inclusion:
$$\xymatrix{{E_2\ }\ar[d]_{\xi}\ar@{^{(}->}[r]^i & X \ar[d]^{\ph_1} \\
T\ar[r] & {Y_1}
}$$
It is easy to see that $i_*(\ker \xi_*)=\ker(\ph_1)_*=\R R_1$. Since
in general $i_*$ is not injective,
  $\xi$ does not need to be an elementary contraction;
however it is birational with
$\Exc(\xi)=E_1\cap E_2$, and $\xi(\Exc(\xi))\subset T$ is a curve.

Notice also that $\rho_T$ is the codimension of $\ker \xi_*$ in
$\N(E_2)$, and since $i_*(\ker \xi_*)$ has codimension $2$ in
$\N(E_2,X)$, 
 we see that $\rho_T\geq 2$,
and $\rho_T= 2$ if and only if $\ker\xi_*\supseteq\ker i_*$.

The diagram 
$$\xymatrix{
{E_2}\ar[r]^{\xi}\ar[d]_{\ph_{2|E_2}}& T\\
A &
}$$
gives a proper, covering family of irreducible rational curves in
$T$, see \cite[\S 5.4]{debarreUT} and references therein. 
This family of curves induces an equivalence relation on $T$ as a set
($E_2$-equivalence in the terminology of  \cite{debarreUT}), where two
points $t_1,t_2\in T$ are equivalent if there exist $F_1,\dotsc,F_m$
fibers of $\ph_{2|E_2}$ such that $\xi(F_1\cup\cdots\cup F_m)$ is
connected and contains both $t_1$ and $t_2$.

By \cite[Th.~5.9]{debarreUT}
there exists a dense open subset $T_0\subseteq T$, closed
for the equivalence relation, and a proper
morphism $\alpha_0\colon
T_0 \to C_0$, where $C_0$ is a normal variety,
such that every fiber of $\alpha_0$ is an equivalence
class.

 Let $S\subset T_0$ be a fiber of $\alpha_0$. We know that
 $\dim\N(S,T)=1$ by \cite[Prop.~IV.3.13.3]{kollar}, and since
 $\rho_T>1$ we know that $S\subsetneq T$ and $\dim C_0>0$.
Moreover 
$\xi^{-1}(S)$ is a union of fibers of $\ph_2$, thus it intersects
$E_1\cap E_2=\Exc(\xi)$, 
so that $S\cap\xi(\Exc(\xi))\neq \emptyset$. Hence $\xi(\Exc(\xi))$
intersects every fiber of $\alpha_0$, which means that $\dim C_0=1$
and every fiber of $\alpha_0$ has codimension $1$. 

Now
if $C$ is the smooth projective curve
containing $C_0$ as an open subset, it is not difficult to see that
the rational map $\alpha_0\colon T\dasharrow C$ extends to a contraction
 $\alpha\colon T\to C$, whose fibers are equivalence
classes, and we get a diagram:
$$\xymatrix{
{E_2}\ar[r]^{\xi}\ar[d]_{\ph_{2|E_2}}&T\ar[d]^{\alpha}\\
A\ar[r]&C
}$$
We deduce that
$\rho_T=2$, and $\ker\xi_*\supseteq\ker i_*$. We refer the interested
reader to \cite{unsplit} and \cite[\S 4]{fanos} for related results.

\medskip

We have 
$$\R(R_1+R_2)\subseteq i_*\bigl(\ker(\alpha\circ\xi)_*\bigr),$$
and since $\ker(\alpha\circ\xi)_*$ is a hyperplane in $\N(E_2)$
and contains $\ker i_*$, its image under $i_*$ must be
$\R(R_1+R_2)$. In particular we see that $\NE((\ph_3)_{|E_2})$ can not
be contained in  $\ker(\alpha\circ\xi)_*$. 

\medskip

Let's show that $E_1\cdot R_3=0$. In fact if $E_1\cdot R_3>0$, then
reasoning as for $R_2$ we get a second contraction $\alpha'\colon T\to
C'$, where $C'$ is another smooth curve. 
Moreover  $\NE((\ph_3)_{|E_2})$ is
contained in $\ker(\alpha'\circ\xi)_*$, hence
$\alpha\circ\xi\neq\alpha'\circ\xi$ and
$\alpha\neq\alpha'$.
However $\dim T=n-1\geq 3$,
and the fibers of $\alpha$ and $\alpha'$ are Cartier divisors which
should intersect only in finitely many points, which is impossible.

Thus $E_1\cdot R_3=0$, $\N(E_1,X)=\R(R_1+R_3)$, and
\stepcounter{thm}
\begin{equation}\label{passi}
\N(E_1,X)\cap\NE(X)=R_1+R_3.
\end{equation}
Then $\N(E_1,X)$ can not contain other extremal rays, and
$R_1+R_3$ is a face of $\NE(X)$ by the following Remark.
\begin{remark}\label{pranzo}
Let $X$ be as in $(\ref{setup}.1)$ and $S_1$ a divisorial extremal ray
of $\NE(X)$
with
exceptional divisor $G_1$,
such that $G_1\cdot S\geq 0$ for every extremal ray $S\neq S_1$. Let
$S_2$ be a birational
 extremal ray of $\NE(X)$ with $G_1\cdot S_2=0$. Then $S_1+S_2$ is a
face of $\NE(X)$, whose contraction is birational. This is probably
well-known; similar properties can be found in \cite{nikulin}. 

Indeed let
$C_i\subset X$ be a curve with $[C_i]\in S_i$ for $i=1,2$. If 
$S_1+S_2$ were not a face of $\NE(X)$, we should have
$$\lambda_1C_1+\lambda_2C_2\equiv\sum_{j=3}^m\lambda_jC_j$$
where $\lambda_j\in\Q_{\,>0}$ for every $j=1,\dotsc,m$, and for $j\geq 3$
$[C_j]$ belongs to an extremal ray $S_j$ with $G_1\cdot S_j\geq
0$. Then intersecting with $G_1$ we get a contradiction.

Moreover if $C\subset X$ is an irreducible curve with
$[C]\in S_1+S_2$, then either $C\cdot G_1<0$, or $[C]\in S_2$, so
that $C\subset G_1\cup \Lo(S_2)$.
\end{remark}

We go on with the proof of \eqref{vernavola}, and
 consider the $3$-dimensional
cone $$i_*(\overline{\NE}(E_2))\subseteq \N(E_2,X)\cap\NE(X),$$
which contains $R_1$, $R_2$, and $R_3$.
Since $R_1+R_3$ and $R_2+R_3$ are faces of $\NE(X)$, they are faces of
$i_*(\overline{\NE}(E_2))$ too. 
On the other hand  $\NE(\alpha\circ\xi)$ is a face of
$\overline{\NE}(E_2)$, and since $\ker(\alpha\circ\xi)_*\supseteq
\ker i_*$, $i_*(\NE(\alpha\circ\xi)$ is a face $S$ of
$i_*(\overline{\NE}(E_2))$, contained in
$i_*(\ker(\alpha\circ\xi)_*)=\R(R_1+R_2)$, and containing both $R_1$
and $R_2$. Therefore $S=R_1+R_2$, and hence
$$i_*(\overline{\NE}(E_2))=R_1+R_2+R_3,$$
which implies \eqref{vernavola}.

\medskip

Now let's consider $\ph_1\colon X\to Y_1$ and the divisor
$\ph_1(E_2)\subset Y_1$. Let $\eta\colon Y_1\to W$
be an elementary contraction
with
$\ph_1(E_2)\cdot \NE(\eta)>0$, as in Rem.~\ref{basic}. 
Moreover let $R_4$ be the extremal ray of $\NE(X)$ such that $R_1+R_4$
is a face and $(\ph_1)_*(R_4)=\NE(\eta)$.

Since $\dim\N(\ph_1(E_2),Y_1)=2$, if $\eta$ is of fiber type we get
$\rho_W\leq 2$ and $\rho_X\leq 4$.

\medskip

Suppose that $\eta$ is birational. Let's show that $\eta$ must be
finite on $\ph_1(E_2)$. If not, there should be curve $C\subset
E_1\cup E_2$ with $[C]\in R_4$. But $[C]\in R_1+R_2+R_3$
by \eqref{vernavola} and \eqref{passi}, which yields
either $R_4=R_2$ or $R_4=R_3$. In both cases we would get
$\Exc(\eta)=\ph_1(E_2)$ and $\ph_1(E_2)\cdot \NE(\eta)<0$, a contradiction.

Thus $\eta$ is finite on $\ph_1(E_2)\supset\ph_1(E_1)$ and must have
fibers of dimension at most~$1$. Then by Lemma
\ref{3/4} $\eta$ is a divisorial Mori
contraction with $\Exc(\eta)\cap \ph_1(E_1)=\emptyset$,
$R_4$ is of type $(n-1,n-2)^{sm}$, and
$E_4\cap E_1=\emptyset$. Moreover $\Exc(\eta)$ must intersect
$\ph_1(E_2)$, so that $E_4\cap E_2\neq\emptyset$.

\medskip

\begin{center}
\scalebox{0.50}{\includegraphics{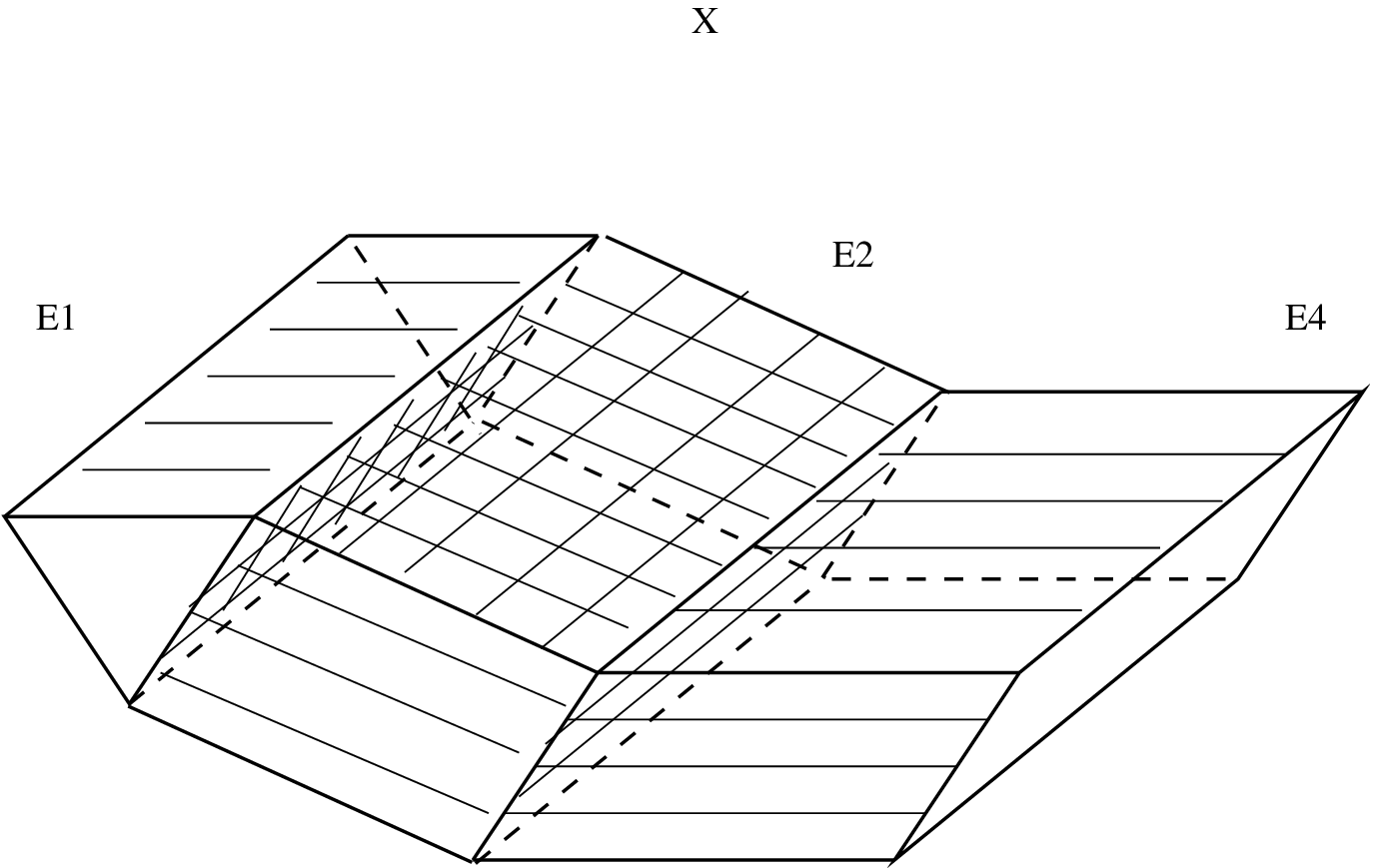}}

\medskip

{\footnotesize Figure 4: the case $\eta$ birational.}

\end{center}

\medskip

Since $E_2$ can not contain curves in $R_4$, we have $E_2\cdot
R_4>0$. 

If $R_2+R_4$ is a face of $\NE(X)$, 
then $(\ph_2)_*(R_4)$ is an extremal ray of $\NE(Y_2)$ with
$D\cdot(\ph_2)_*(R_4)>0$, whose locus is either $\ph_2(E_4)$ or the
whole $Y_2$. However if $C\subset X$ is a 
non trivial fiber of $\ph_4$, it is easy
to see that $\ph_2(C)\cdot \ph_2(E_4)\geq 0$, thus the contraction of 
 $(\ph_2)_*(R_4)$ is of fiber type and as before we get 
$\rho_{Y_2}\leq 3$ and $\rho_X\leq 4$.

\medskip

Finally let's assume that $R_2+R_4$ is not a face of $\NE(X)$, and
consider the divisor $\ph_4(E_1)\subset Y_4$.
There exists an extremal
ray $S$ of $\NE(Y_4)$ with $\ph_4(E_1)\cdot S>0$. Let $R_5$ be the
extremal ray of $\NE(X)$ such that $R_4+R_5$ is a face of $\NE(X)$ and 
$(\ph_4)_*(R_5)=S$. We observe that by construction $R_5\neq R_2$.
Since
$\ph_4^{-1}(\ph_4(E_1))=E_1$, we have $E_1\cdot R_5>0$, hence $R_5\neq
R_1$ and $R_5\neq R_3$.

Now we apply what we proved so far
to $R_5$. Notice that
$R_5\not\subset\N(E_1,X)$, in particular $R_5$ can not be small.
Then either $\rho_X\leq 4$, or $R_5$ is of type $(n-1,n-2)^{sm}$ with
$E_5\cdot R_1>0$, and there exists a
divisorial extremal ray $R_6\neq R_5$  such that $E_5\cdot R_6<0$ 
and $E_1\cdot R_6=0$. We show that this last case is impossible.

In fact we have
$R_6\subset\N(E_1,X)=\R(R_1+R_3)$ and $R_6\neq R_1$ because they are
of 
different types, so
the only possibility is that $R_6=R_3$
and $E_5=E_2$. 
If $C$ is a curve with numerical class in $R_5$, then $C\subset E_2$,
hence
$[C]\in R_1+R_2+R_3$ by \eqref{vernavola}. 
But $R_5$ is distinct from $R_1,R_2,R_3$,
so we get a contradiction.

\medskip

We still have to consider the case where  $\ph_2$ is of type
$(n-2,n-4)$, and
there exists a smooth prime divisor
$D'\subset X$, disjoint from $E_2$ and having
 a $\pr^1$-bundle structure $\xi\colon D'\to W$, such that
for any fiber $f$ of $\xi$
we have  $D'\cdot f=-1$ and $E_1\cdot f>0$.

Notice that every non trivial fiber of $\ph_1$ must intersect
$E_2$, thus it can not be contained in $D'$. This implies that
$D'\cdot R_1>0$, so that $D'$ intersects every curve contracted by
$\ph_1$. Again since $D\cap E_2=\emptyset$, we see that 
  $\ph_1$ is finite on $E_2$.
This gives
$$n-3=\dim (E_1\cap E_2)=\dim\ph_1(E_1\cap E_2)\leq 1,$$
hence $n=4$. 

\medskip

\begin{center}
\scalebox{0.50}{\includegraphics{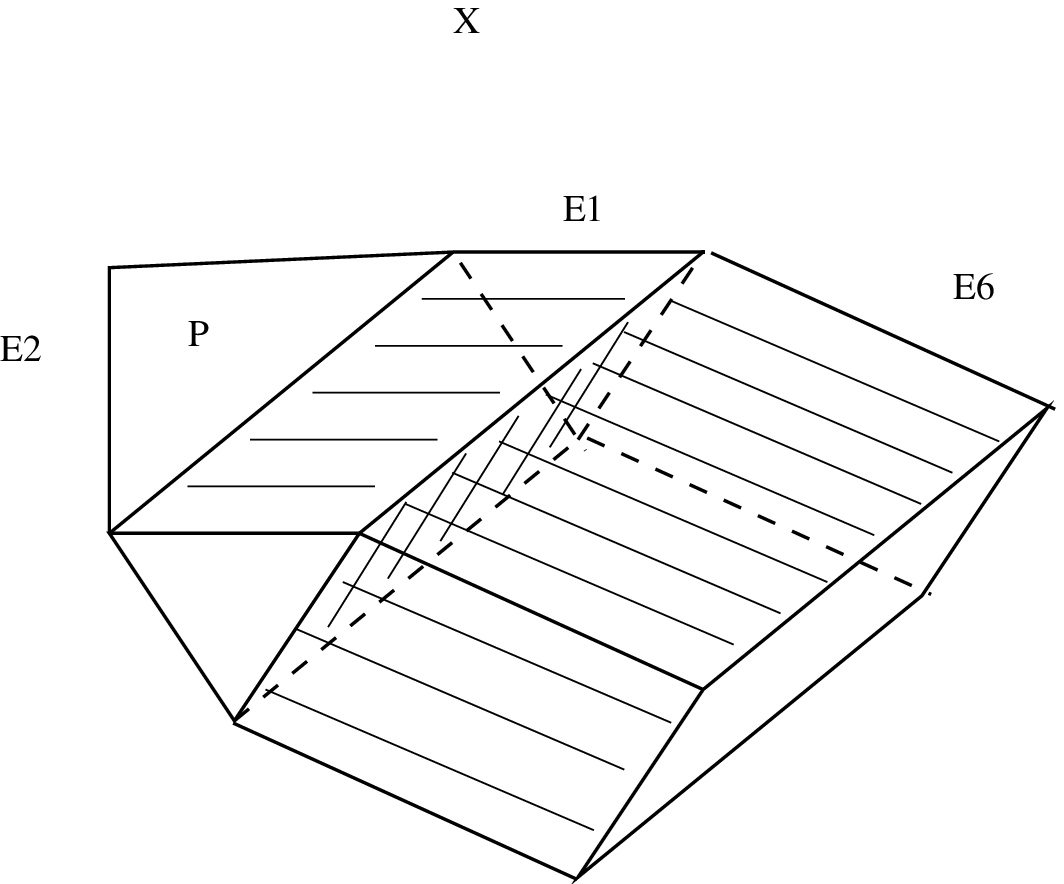}}

\medskip

{\footnotesize Figure 5: $\ph_2$ is $(2,0)$, $n=4$.}

\end{center}

\medskip

We have $\N(E_1,X)\cap\NE(X)=R_1+R_2$, and
$D'\cdot R_1>0$, $D'\cdot R_2=0$, $D'\cdot f<0$.
Thus $[f]\not\in\N(E_1,X)$ and $\rho_X\geq 3$.
Moreover $\xi(D'\cap E_1)=W$, hence
$$\N(D',X)=\R[f]\oplus \N(D'\cap E_1,X)=\R[f]\oplus\R R_1\oplus\R R_2$$
and $\dim\N(D',X)=3$.

Since $D'\cdot f<0$, there exists some extremal ray $\widetilde{R}_2$
of 
$\NE(X)$
with $D'\cdot \widetilde{R}_2<0$. If $\widetilde{R}_2$ 
were small, by \cite{kawsmall}
its exceptional locus would contain $F\cong\pr^2$. Then $\xi(F)=W$,
which would give  $\dim\N(D',X)=2$, a contradiction. Thus
$\widetilde{R}_2$ is
divisorial, with exceptional divisor $D'$. Since
$\widetilde{R}_2\not\subset\N(E_1,X)$, $\widetilde{R}_2$ 
is of type $(3,2)^{sm}$, and
$E_1\cdot \widetilde{R}_2>0$. Applying to
$\widetilde{R}_2$ what we have already proved we get $\rho_X\leq 4$. In fact it
is not difficult to see that $\widetilde{R}_2$ contains $[f]$. 
\end{proof}
\begin{remark}\label{vaccinazione}
Let $X$ and $R_1$ be as in $(\ref{setup}.1)$, 
and suppose that $R_2$ is a birational
 extremal ray with $E_1\cdot R_2>0$, $E_2\cdot
R_1=0$, and $R_2\not\subset\N(E_1,X)$ (see figure 1 on p.~\pageref{fig1}).

We have seen
 in the proof of Th.~\ref{positive} that $R_2$ is of
type $(n-1,n-2)^{sm}$ and $\N(E_2,X)=\R(R_1+R_2)$.

Therefore $E_2\cdot S\geq 0$ for
every extremal ray $S\neq R_2$, and $Y_2$ is Fano by
\cite[Prop.\ 3.4]{wisn}. Moreover by Rem.~\ref{pranzo}
$R_1+R_2$ is a face 
of $\NE(X)$, whose contraction is birational. Notice that the
contraction of $R_1+R_2$
can not send $E_1$ to a point, otherwise we would have 
$\N(E_1,X)=\R(R_1+R_2)$ which is excluded by our assumptions. Thus
$(\ph_2)_*(R_1)$ is an extremal ray of $\NE(Y_2)$, whose contraction
is birational and can not send $\ph_2(E_1)$ to a point. This means
that $Y_2$ has an elementary contraction of type $(n-1,1)$ given by
$(\ph_2)_*(R_1)$, with exceptional divisor $\ph_2(E_1)$,
 and
$\ph_2$ is the blow-up of a smooth fiber of such contraction.
\end{remark}
\begin{proposition}\label{11/4}
Let $X$  and $R_1$ be as in $(\ref{setup}.1)$, and 
suppose that there exists an extremal
ray $R_0\neq R_1$ with $E_1\cdot R_0<0$. 

Then
$\rho_X\leq 5$, $R_0+R_1$ is a face of $\NE(X)$,
$E_1\cong W\times\pr^1$ where $W$ is smooth and Fano, 
$Y_0$ is smooth, and
$\ph_0$ is the blow-up of a smooth subvariety isomorphic to $W$.

If moreover $\rho_X=5$, then there exists a smooth Fano variety $Z$
with $\rho_Z=3$ and $\dim Z=n$, having an elementary contraction of
type $(n-1,1)$, 
 such that $X$ is the blow-up of $Z$ in two fibers of such
 contraction.
\end{proposition}
\begin{proof}
Every non trivial
fiber of $\ph_0$ is contained in $E_1$ and hence has
dimension $1$. Therefore $R_0$ is of type $(n-1,n-2)^{sm}$,
$Y_0$ and $E_1$ are smooth, $\ph_0$ is the blow-up of a smooth,
codimension $2$ subvariety $W\subset Y_0$, and $E_1$ is
a $\pr^1$-bundle over $W$.

Moreover $\N(E_1,X)=\R(R_0+R_1)$, $\N(E_1,X)\cap\NE(X)=R_0+R_1$,
and there are no other extremal
rays with negative intersection with $E_1$.

For $i=1,2$ let $C_i$ be a curve in $R_i$ and $H_i$ a nef divisor such
that for every extremal ray $S$ of $\NE(X)$, $H_i\cdot S=0$ if and
only if $S=R_i$.
The divisor
$$H:=(H_0\cdot C_1)H_1+(H_1\cdot C_0)(-E_1\cdot C_1)H_0+(H_0\cdot
C_1)(H_1\cdot C_0)E_1$$
is nef, and for every extremal ray $S$ of $\NE(X)$, $H_i\cdot S=0$ if and
only if $S=R_0$ or $S=R_1$. Thus $R_0+R_1$ is a face of $\NE(X)$.

\medskip

Let's show that $E_1$ is Fano. If $\gamma\in\overline{\NE}(E_1)$ is
non zero, then
$$-K_{E_1}\cdot\gamma=-(K_X+E_1)_{|E_1}\cdot\gamma=-(K_X+E_1)\cdot
i_*(\gamma),$$
where $i\colon E_1\hookrightarrow X$ is the 
inclusion. First of all we observe that $i_*(\gamma)$ is non
zero. Indeed if $A$ is an ample divisor on $X$, then
$$A\cdot i_*(\gamma)=A_{|E_1}\cdot\gamma>0.$$
Moreover
$i_*(\overline{\NE}(E_1))\subseteq\NE(X)$, so that $i_*(\gamma)\in
R_0+R_1$ and hence $E_1\cdot  i_*(\gamma)<0$. This gives
$-K_{E_1}\cdot\gamma>0$. 

\medskip

The restriction $\ph_{1|E_1}\colon E_1\to \ph_1(E_1)$ is surjective
with connected fibers. Since $\ph_1(E_1)$ is covered by fibers of
$\ph_{0|E_1}$, it is a rational curve, and  $\ph_{1|E_1}$ induces a
Mori contraction
$$\phi\colon E_1\la\pr^1$$
which does not contract the fibers of $\ph_{0|E_1}$. 
Then $E_1\cong W\times\pr^1$ by the following Lemma.
\begin{lemma}\label{XXX}
Let $E$ be a smooth variety and $\pi\colon E\to W$ be a
smooth morphism with fiber $\pr^r$. 
Suppose that $E$ has a Mori contraction $\phi\colon
E\to\pr^r$ which is finite on fibers of $\pi$. Then $E
\cong W\times\pr^r$. 
\end{lemma}
We postpone the proof of Lemma \ref{XXX} and carry on with the proof
of Prop.~\ref{11/4}.  

Let $R_2$ be an extremal ray of $\NE(X)$ with $E_1\cdot R_2>0$. Then
$R_2$ is different from $R_0$ and $R_1$, so that
$R_2\not\subset\N(E_1,X)$, and $\ph_2$ is finite on $E_1$ (notice that
necessarily $\rho_X\geq 3$).

If $\ph_2$ is of fiber type, then it is a conic bundle,
$\rho_{Y_2}= 2$ and $\rho_X=3$.

\medskip

Suppose that $\ph_2$ is birational. Then it is of type $(n-1,n-2)^{sm}$,
so $Y_2$ is smooth and $\ph_2$ is the blow-up of
$A:=\ph_2(E_2)\subset\ph_2(E_1)\subset Y_2$. We set $D:=\ph_2(E_1)$.

Notice that $\ph_2(E_1\cap
E_2)=\ph_2(E_2)$, and $C\cdot E_2\geq 0$ for every curve $C\subset
E_1$. Since $\ph_2^*(-K_{Y_2})=-K_X+E_2$, using the projection formula
we see that $Y_2$ is Fano.

Let $\psi\colon Y_2\to Z$ be an elementary contraction such that 
$D\cdot
\NE(\psi)>0$, as in Rem.~\ref{basic}.
If $\psi$ is
of fiber type, then $\rho_Z\leq 2$ and $\rho_X\leq 4$.

\medskip

Assume that $\psi$ is birational. 
Then $\psi$ must be finite on
$D$, because
\begin{equation} \label{seda}\begin{split}
\N(D,Y_2)\cap\NE(Y_2)=&(\ph_2)_*\bigl(\N(E_1,X)\cap\NE(X)  \bigr)\\=
&(\ph_2)_*(R_0)+(\ph_2)_*(R_1).\end{split}\end{equation}
If $\psi$ were not finite on $D$, it should 
be $\NE(\psi)=(\ph_2)_*(R_0)$
or $\NE(\psi)=(\ph_2)_*(R_1)$; in both cases $\Exc(\psi)=D$, which
contradicts $D\cdot
\NE(\psi)>0$.
Thus $Z$ is smooth and $\psi$ is of type $(n-1,n-2)^{sm}$. 

Lemma \ref{3/4} says that $\Exc(\psi)\cap
A$ is a union of fibers of $\psi$, but $\psi$ is finite on $A$, so 
$\Exc(\psi)\cap
A=\emptyset$. 
Hence the composition
$$\psi\circ\ph_2\colon X\la Z$$
is just the blow-up of two disjoint subvarieties in $Z$.
Set $\widetilde{E}_2:=\ph_2^{-1}(\Exc(\psi))$, so that
$\Exc(\psi\circ\ph_2)=E_2\cup \widetilde{E}_2$.

Let's show that $E_2\cdot R_1=\widetilde{E}_2\cdot R_1=0$.
The intersection $E_1\cap {E}_2$ has pure dimension $n-2\geq 2$, thus 
$\ph_{1|
E_1\cap {E}_2}\colon E_1\cap {E}_2\to\ph_1(E_1)$ has positive dimensional
fibers. Take a curve $C$ in one of these fibers: then $[C]\in R_1$ and
$C\subset E_2$, thus $C\cap \widetilde{E}_2=\emptyset$, so
$\widetilde{E}_2\cdot R_1=0$. In the same way we see that $E_2\cdot
R_1=0$.

Therefore both $E_1\cap E_2$ and $E_1\cap\widetilde{E}_2$ are union of
finitely many fibers of $\ph_1$. 

\medskip

\begin{center}
\scalebox{0.50}{\includegraphics{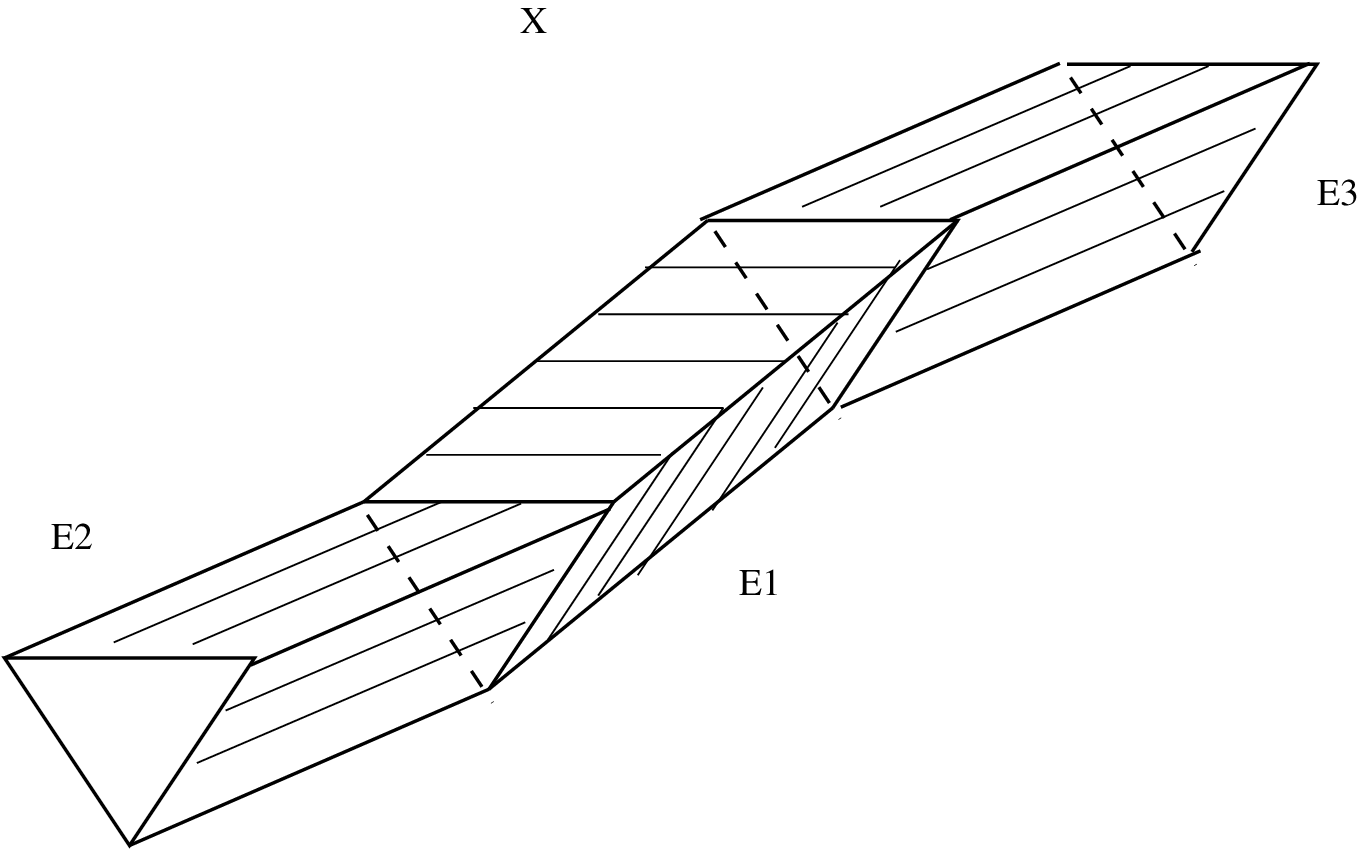}}

\medskip

{\footnotesize Figure 6: the case $\psi$ birational.}
\end{center}

\medskip

We apply Rem.~\ref{vaccinazione} to $R_1$ and $R_2$, and deduce
that $R_1+R_2$ is a face of $\NE(X)$, and
$S_1:=(\ph_2)_*(R_1)$ is an extremal ray of $\NE(Y_2)$ of type
$(n-1,1)$, with exceptional divisor $D$.

By \eqref{seda} we know that apart from $S_1$, the other possible
extremal ray contained in $\N(D,Y_2)$ is $(\ph_2)_*(R_0)$. 
It is easy to see that
$E_2\cdot R_0>0$ and $D\cdot (\ph_2)_*(R_0)\geq 0$. 
This shows
that $S_1$ is the unique extremal
ray of $\NE(Y_2)$ having negative intersection with $D$, and
Th.~\ref{positive} yields $\rho_{Y_2}\leq 4$ and
$\rho_X\leq 5$.

Recall that $\NE(\psi)$ is a
birational extremal ray of $\NE(Y_2)$ with
$D\cdot \NE(\psi)>0$ and 
$\NE(\psi)\not\subset \N(D,Y_2)$. 
Moreover
$\widetilde{E}_2\cdot R_1=0$ in $X$ yields
$\Exc(\psi)\cdot S_1=0$ in $Y_2$. 
Then we can apply Rem.~\ref{vaccinazione} to $Y_2$, $S_1$, $\NE(\psi)$
as we did for $X$, $R_1$, $R_2$. We deduce that $Z$ is Fano,
$\psi_*(S_1)$ is an extremal ray of type $(n-1,1)$ with exceptional
divisor $\psi(D)$, and $X$ is the blow-up of $Z$ in two
fibers of the associated 
contraction. Notice that $\psi\circ\ph_2$ is finite and
birational on $E_1$, thus the normalization of $\psi(D)$ is
$W\times\pr^1$.
\end{proof}
\begin{proof}[Proof of Lemma \ref{XXX}]
We proceed similarly to
the proof of \cite[Lemma 1.2.2]{AW2001}.
Let $K_{\pi}$ be the relative canonical bundle of $\pi$. Let's show
that $K_{\pi}\cdot C=0$ for every curve $C$ contracted by $\phi$.
Since $\phi$ is a Mori contraction, it is enough to show this when $C$
is an irreducible rational curve. Then $\pi(C)$ is again an 
irreducible rational curve. Let $\nu\colon\pr^1\to W$ be the morphism
given by the normalization of $\pi(C)\subset W$, and consider 
the fiber product:
$$\xymatrix{{E_C}\ar[r]^{\tilde{\nu}}\ar[d]_{\pi_C} & {E}\ar[d]^{\pi} \\
{\pr^1}\ar[r]^{\nu} &{W}
}$$
where $E_C\to\pr^1$ is a $\pr^r$-bundle. Notice that
$K_{\pi_C}=\tilde{\nu}^*(K_{\pi})$. 

Let $\phi_C$ be the composition given by the following diagram:
$$\xymatrix{
{E_C}\ar[r] \ar@/^1pc/[rrr]^{\phi_C}&
{\pi^{-1}(\pi(C))}\ar[rr]_(.6){\phi_{|\pi^{-1}(\pi(C))}}
&&{\pr^r}
}$$
Then $\phi_C$ is surjective and its Stein factorization gives a
contraction $\xi_C\colon E_C\to P$ which is finite on fibers of
$\pi_C$, and such 
that $\dim P=r$. This easily implies (for instance using toric
geometry) that $E_C\cong\pr^1\times\pr^r$, $P\cong\pr^r$, and $\xi_C$
is the projection. 
Then $K_{\pi_C}=\xi_C^*(K_{\pr^r})$.

Now set $\widetilde{C}:=\tilde{\nu}^{-1}(C)\subset E_C$. Since
$\phi(C)=\{pt\}$, we have $\phi_C(\widetilde{C})=\{pt\}$ and hence
$\xi_C(\widetilde{C})=\{pt\}$. Moreover
$\tilde{\nu}_*(\widetilde{C})=mC$ for some $m\in\Z_{\geq 1}$. Finally
$$K_{\pi}\cdot C=\frac{1}{m}K_{\pi}\cdot \tilde{\nu}_*(\widetilde{C})=
\frac{1}{m} \tilde{\nu}^*(K_{\pi})\cdot \widetilde{C}=
\frac{1}{m}K_{\pi_C}\cdot \widetilde{C}=\frac{1}{m}\xi_C^*(K_{\pr^r})
\cdot \widetilde{C}=0.$$

\medskip

Now let $F$ be a general fiber of $\phi$ and let $d$ be the degree of
the finite map $E\to W\times\pr^r$ induced by $\pi$ and $\phi$. Then 
$g:=\pi_{|F}\colon F\to W$ is finite of degree $d$. Moreover
$$K_F=(K_E)_{|F}=(K_{\pi}+\pi^*K_W)_{|F}=(K_{\pi})_{|F}+g^*K_W.$$
Since $K_{\pi}$ is numerically trivial on $F$ and $F$ is Fano, we have 
$(K_{\pi})_{|F}\cong\mathcal{O}_F$, so that $K_F=g^*K_W$ and $g$ is
\'etale.
Then $W$ is Fano too, in particular it is simply connected, thus $g$
is an isomorphism and $d=1$.
\end{proof}
\begin{proof}[Proof of Cor.~\ref{dim4}]
The statement is a straightforward consequence of Th.~\ref{5} and
\cite[Th.~1.1]{takagi}. 
\end{proof}
\begin{proof}[Proof of Cor.~\ref{4folds}]
Let $X$ be a smooth Fano $4$-fold with $\rho_X\geq 7$. Then $X$ can
not have elementary contractions of type $(3,0)$,
$(3,1)$,  or $(4,2)$ by Th.~\ref{5}, Prop.~\ref{1}, and
\cite[Cor.~1.2]{fanos}. 
Thus the possibilities are just $(3,2)$,
$(2,0)$, or $(4,3)$.
But if $X$ has an elementary contraction of type $(4,3)$, then 
\cite[Cor.~1.2]{fanos} implies that either
$X\cong\pr^1\times\pr^1\times S$, or $X\cong\mathbb{F}_1\times S$,
where $S$ is a Del Pezzo surface; in particular $\rho_X=2+\rho_S\leq
11$. Therefore we have the statement.
\end{proof}
\begin{example}\label{ultimo}
It is not difficult to find examples 
of Fano varieties $X$ as in Th.~\ref{5},
with $\rho_X=5$. For instance in the toric case, we know
after \cite{sato2} (and \cite{bat2} for the $4$-dimensional case) 
that there are exactly $n-2$ possibilities for $X$, which can be
obtained as follows.

Let $a$ be an integer with $1\leq a\leq n-2$ and consider
$$Z:=\pr_{\pr^{n-2}\times\pr^1}(\mathcal{O}(0,1)
\oplus\mathcal{O}(a,0)).$$
Then $Z$ is Fano with $\rho_Z=3$.
The $\pr^1$-bundle $Z\to \pr^{n-2}\times\pr^{1}$ has a section $E_Z$
with normal bundle
$\mathcal{N}_{E_Z/Z}\cong\mathcal{O}_{\pr^{n-2}\times\pr^{1}}
(-a,1)$, and $Z$ has an extremal ray of type $(n-1,1)$ with
exceptional divisor $E_Z$. Blowing-up $Z$ along
$\pr^{n-2}\times\{p_1,p_2\}\subset E_Z$ (where $p_1,p_2\in\pr^1$ are
two distinct points) yields a toric Fano variety $X$ with $\rho_X=5$,
where the proper transform ${E}\cong\pr^{n-2}\times\pr^1$
of $E_Z$ has normal bundle 
$\mathcal{O}_{\pr^{n-2}\times\pr^{1}}
(-a,-1)$. Finally $X$ has an extremal ray of type $(n-1,1)$ and one of
type $(n-1,n-2)^{sm}$, both with exceptional divisor ${E}$.
\end{example}
\footnotesize
\providecommand{\bysame}{\leavevmode\hbox to3em{\hrulefill}\thinspace}
\providecommand{\MR}{\relax\ifhmode\unskip\space\fi MR }
\providecommand{\MRhref}[2]{%
  \href{http://www.ams.org/mathscinet-getitem?mr=#1}{#2}
}
\providecommand{\href}[2]{#2}

\end{document}